# On the Riemann zeta-function, Part I: Outline

By Anthony Csizmazia

E-mail: apcsi2000@yahoo.com

## Abstract

Results of a multipart work are outlined. Use is made therein of the conjunction of the Riemann hypothesis, RH, and hypotheses advanced by the author. Let z(n) be the nth nonreal zero of the Riemann zeta-function with positive imaginary part in order of magnitude thereof. A relation is obtained, of the pair z(n) and the first derivative thereat of the zeta-function, to the preceding such pairs and the values of zeta at the points one-half plus a nonnegative multiple of four. That relation is derived from two forms of the density of the Laplace representation, on a certain vertical strip, of a meromorphic function constructed from zeta. Specific functions which play a central role therein are proven to have analytic extensions to the entire complex plane. It is established that the Laplace density is positive. That positivity implies RH and that each nonreal zero of zeta is simple. A metric geometry expression of the positivity of the density is derived. An analogous context is delineated relative to Dirichlet L-functions and the Ramanujan tau Dirichlet function.

**MSC (Mathematics Subject Classification).** 11Mxx Zeta and L-functions: analytic theory. 11M06 $\zeta(s)$ and $L(s, \chi)$. 11M26 Nonreal zeros of $\zeta(s)$ and $L(s, \chi)$; Riemann and other hypotheses. 11M41 Other Dirichlet series and zeta functions. 30xx Functions of a complex variable. 44A10 Laplace transform. 42xx Fourier analysis. 42A38, 42B10 Fourier and Fourier-Stieltjes transforms and other transforms of Fourier type. 42A82 Positive definite functions. 60E10 Characteristic functions; other transforms. 33C15 Confluent hypergeometric functions, Whittaker functions, $_1F_1$.

**Journal of Number Theory classifications.** 120.000 Zeta and L-functions 120.010 Analytic theory 120.020 Zeros of L-functions 120.030 Riemann and other hypotheses

# On the Riemann zeta-function

## Table of contents.

**Part I: Outline**

**Abstract. Keywords. MSC (Mathematics Subject Classification). Journal of Number Theory classifications.**

**Index of abbreviations. Index of symbols.**

**Introduction**

**§1 Definitions**

The functional equation of $\zeta(s)$. Symmetries of $\zeta(s)$, $\xi(s)$, $n(s)$, $f(s)$. Zeros of $\zeta(s)$, $\xi(s)$, $n(s)$. Poles of $f(s)$. The nonreal zeros of $\zeta(s)$.

**§2 Standard conjectures**

(2.1) The Riemann Hypothesis, RH
(2.2) The simple zeros conjecture, SZC
(2.3) The Lindelöf hypothesis, LH
(2.4) The Lehmer phenomenon

**§3 Context and statement of the main unconditional theorem**

Integrability of $|f(s)|$ on vertical lines outside of the critical strip. Definition and properties of $F(s)$. Main unconditional theorem. The alteration of the density under a transition of domain of a Laplace representation. A geometric consequence of the Main unconditional theorem (4).

**§4 Conjectures introduced by the author**

(4.k) Conjecture k, k = 1, 2, 3, 4, 5
(4.1) The asymptotic behavior of $|b(s)|$ on a vertical strip of finite width.
(4.2) Note on the finiteness of A.

**§5 Context and statement of the main conditional theorem**

Introduction
(5.1) Analyticity of $f(s)$
(5.2) Integrability of $|f(s)|$ on vertical lines

## §6 Relations of $\gamma_n$, $\zeta'(\frac{1}{2} + i\gamma_n)$ to their predecessors and the $\zeta(\frac{1}{2} + 4k)$.

A key observation. Relations of $\gamma_n$, $\zeta'(\frac{1}{2} + i\gamma_n)$ to their predecessors and the $\zeta(\frac{1}{2} + 4k)$. The Hadamard factorization of $\xi(\frac{1}{2} + s)$. $\xi$-fission. Laplace transform representation of $ip_{i,+}(is)$. Representation of $p_{i,+}(z)$ via j(y).

## §7 f(s) as a meromorphic characteristic function on C

The Main conditional theorems (2)-(3). The Main conditional theorem (2) implies RH and SZC. The role of b(s) in the groove property of n(s). Polynomial counterexamples for the groove property. Metric norms and analytic characteristic functions.

## §8 Meromorphic characteristic functions arising from $L(s, \chi)$ or r(s)

(8.1) Dirichlet L-functions, $L(s, \chi)$
(8.2) The Ramanujan tau Dirichlet function, r(s)

## Appendix

Kronecker's approximation theorem and Dirichlet L-functions, $L(s, \chi)$.

## Acknowledgements

## References

## Part II: The Laplace representation of $1/(\sin(\pi s/4)\cdot 2\xi(\frac{1}{2} + s))$.

## Abstract. Keywords. MSC (Mathematics Subject Classification). Journal of Number Theory classifications.

## Review of elements of Part I.

## §1 The role of $\sin(\pi s/4)$ in $f(s) := 1/(\sin(\pi s/4)2\xi(\frac{1}{2} + s))$.

Properties of $\zeta(s)$. Property of $\zeta(s)$. Integrability of $|m(x + it, \beta)|$ in t. The Laplace representation of $m(z, \beta)$. The Laplace representation of $\pi/\sin(\pi s)$. Translation principal.

**§2 The Laplace representation on $V_0'$ of $f(s) := 1/(b(s)\zeta(\frac{1}{2} + s))$ via that of $1/b(s)$.**

**§3 Strategy for determining the Laplace representation on $V_0'$ of $1/b(s)$.**

Divergence of the formal partial fraction expansion of $1/b(s)$. Integrability of $|f(s)|$ on vertical lines outside of the critical strip. Bound for $1 / |\zeta(z)|$. Integrability of $|N(z, \frac{1}{4})|$ on vertical lines. Determination of the Laplace representations of $f_0(s)$ and $1/b(s)$ on $V_0'$ from that of $F(z, \frac{1}{4})$ on $V(\frac{1}{4}, 2)$. Translation relation for $F(z, \beta)$. Translation relation for $f(z)$. The Mellin transform representation of $j(u, m)$. The Mellin transform representation of $1/(z)_m$. The Mellin transform representations of $1/(S + 2m - \frac{1}{2})$, $j(\frac{1}{2} + S, m)$ and $f_0(S)$. The Mellin transform representation of $f(s)$, for s on $V_0' + 4w$. Positivity. Determination of the Laplace representation of $f(s)$ on $V_{4w}'$. Determination of the Laplace representation of $F(z, \beta)$ relative to z on $V_{4w}$. Integrability of $|F(z, \beta)|$ on vertical lines. Divergence of the formal partial fraction expansion of $F(z, \beta)$. Splitting $F(u, \beta, 1)$. Trigonometric identity. The splitting of $F(u, \beta, 1)$ via $E(u, \beta, 1)$ and $E(u, \beta, 2)$. Trigonometric asymptotics. The asymptotic behavior of $|E(u, \beta, 1)|$ and of $|E(u, \beta, 2)|$. Integrability. The Mellin transform representation of $E(u, \beta, 2)$. The Mellin transform representation of $E(u, \beta, 1)$. Definition and properties of $W(z, \beta)$. Translation relation for $W(z, \beta)$. Mellin transform representation of $((\pi/2)/\cos((\pi/2)(u + \beta)))\cdot\Gamma(u)$. Definition and properties of $B_0(z, \beta)$. Mellin transform representation of $E(u, \beta, 1)$. The Mellin transform representation of $F(u, \beta, 1)$. Definition and properties of $M(z, \beta)$. Mellin transform representation of $F(u, \beta, 1)$. Definition and properties of $I(p, z, u)$. The determination of $W(z, 1 + \beta)$ from $I(p, z, u)$. The determination of $I(p, z, u)$ from $I(p, z/u)$. The determination of $W(z, 1 + \beta)$ from $I(-\beta, \pm iz)$.

**§4 The incomplete gamma functions and confluent hypergeometric functions**

Relation between $I(p, z)$ and $\Gamma(1 - p, z)$. The confluent hypergeometric functions $M(a, B, z)$ and $U(a, a, z)$. Kummer's equation. Relation between $\varphi(1 + \beta, z)$ and $\gamma(\beta, z, *)$. The Laplace representation of $\varphi(1 + \beta, z)$. Relation between $\varphi(1 + \beta, z)$ and $\varphi(\beta, z)$. Relation between $I(-\beta, z)$ and $\varphi(1 + \beta, z)$.

**§5 $H(z, \beta)$ and the Mellin transform representation of $F(z, \beta)$.**

The determination of $W(z, 1 + \beta)$ from $H(z, \beta)$. The determination of $M(z, \beta)$ from $W(z, 1 + \beta)$. The determination of $M(z, \beta)$ from $H(z, \beta)$. Integrability of $|F(z, \beta)|$ on vertical lines. Representation of $H(z, \beta)$. Uniform boundedness. Positivity. Translation relation for $H(z, \beta)$. The determination of $H(z, \beta)$ from $W(z, \beta - 1)$. The representation of $W(z, \beta)$ via $R(z, \beta)$. The boundedness of $H(z,$ β). Mellin transform representation of F(z, β). Positivity.

## §6 The Mellin transform representation of f(s, β) := 1/(sin(πs/4)·2ξ(2β + s)).

Relations used to determine the Mellin transform representation of f(s, β). Integrability of 1/|b(s, β)| on vertical lines. The Mellin transform representation of $f_0(s, β)$. The Mellin transform representation of 1/b(s, β). The Mellin transform representation of f(s, β). The determination of $P_{4w}$ from $T_0$. The order of $P_0(r, β)$ in r. Results when β = ¼ . Main unconditional theorem (1), (4) (i). Positivity when β and w are nonnegative. Results when β = ¼ . Main unconditional theorem (4) (i), (ii). Metric norms and analytic characteristic functions. Metric result when β = ¼ .

## References

## Part III: The partial fraction representation of 1/(sin(πs/4)2ξ(½ + s)).

## Abstract. Keywords. MSC (Mathematics Subject Classification). Journal of Number Theory classifications.

## Review of elements of Parts I-II, IV.

## Introduction

The conditional implications (*), (**). A general lemma on partial fraction expansions: Absolute convergence, Vanishing at infinity.

## §1 Unconditional results.

Assumptions and definitions. Bounds on $|ζ(z)|^{-1}$. A bound for the remainder in a Taylor series. The assumptions for Claims 1.1, 1.2, 2.3 and 2.4. The special cases (′) and (″).

## §2 Conditional results.

Counterexample. The Hadamard factorization of ξ(½ + s). The monotonicity principle. Conditional theorem 2.1: Vertical vanishing of f(s) − T(s) on $B_i(α)$. Conditional theorem 2.2: Partial fraction representation of f(s).

## **Appendix** Complete monotonicity.

## References

**Part IV: On the Riemann zeta-function and meromorphic characteristic functions.**

**Abstract. Keywords. MSC (Mathematics Subject Classification). Journal of Number Theory classifications.**

**Review of elements of Parts I-III.**

**§1 The conditional Laplace representation of f(s) on $V_0$.**

**§2 Proof of the Main conditional theorem (1).**
Main conditional theorem (1)
(i) The equality of the conditional and unconditional Laplace densities.
(i') The boundedness of the density.
(ii) The conditional extension of the unconditional Laplace representation of f(s) on $V_0'$ to $V_0$.

**§3 Proofs of the Main conditional theorems (2)-(3).**

**§4 Metric norms and analytic characteristic functions.**

**References**

**Part V: On the Riemann zeta-function: A relation of its nonreal zeros and first derivatives thereat to its values on ½ + 4N.**

**Abstract. Keywords. MSC (Mathematics Subject Classification). Journal of Number Theory classifications.**

**Review of elements of Parts I, IV.**

**§1 Proof of the conditional relations of $\gamma_n$, $\zeta'(½ + i\gamma_n)$ to their predecessors and the $\zeta(½ + 4k)$.**

**§2 Representation of $p_{i,+}(z)$ via j(y).**

**References**

**Part VI:** A generalization of the Fourier metric geometries of J. von Neumann and I. J. Schoenberg.

**Abstract. Keywords. MSC (Mathematics Subject Classification). Journal of**

**Number Theory classifications.**

**Review of elements of Parts I, II and IV.**

**Introduction**

**§1**

**§2**

**References**

**Index of abbreviations**

**Part I**

**§2**
**(2.1) (**Riemann Hypothesis**)** RH.
**(2.2) (**Simple zeros conjecture**)** SZC.
**(2.3)** (Lindelöf hypothesis) LH.
**(2.4)** (Gaussian unitary ensemble) GUE.

**§8**
**(8.1)** (Dirichlet L-function conjecture) DLFC. (Simple zeros conjecture) SZCD.
**(8.2)** Ramanujan conjecture, RC. Simple zeros conjecture, SZCR.

**Appendix**

d-invasive.

**Part II**

**§1**
**Corollary 1.1** (Analytic characteristic function) ACF.

**Index of symbols**

(Complex plane) C. (Real line) **R**, $\mathbb{R}$.

**Part I**

**§1**
**(1.1)** (Functions) l(s), a(s), ξ(s), n(s), f(s), b(s). (Vertical strips) V(x₀, x₁), V[x₀,

$x_1$], V($\varepsilon$).

**§2**
**(2.2)** $\gamma_n$
**(2.4)** $\delta_k$

**§3**
(Vertical strips) $V_u$, $V_u'$. (Pochhammer symbol) $(z)_n$. (Coefficients) $\tilde{c}(4k)$, $c(z)$. (Entire function) $P_0(z)$. (Open disk) $B(z, r)$. (Positive integers) **N**. (partial fraction expansions) $F(s)$, $p_r(s)$. (Integers) **Z**. (Sets of zeros) $Z_i$, $Z°$. **(**Laplace density function) $g(z)$. (Entire function in z) $P_{4w}(z)$. (Metric norm) $m_x(t)$. (Metric) $d_x(t_1, t_2)$.

**§4**
(Conjecture k) Ck. (Conditional Laplace density) $g_0(y)$.
**(4.1)** Definitions $\delta_k'$. (Open interval) $I_k(\alpha)$. (Semicircle) $S_k(\alpha)$. (Boundary) $T_*(\alpha)$. $t(1)$ $x(t, \alpha)$ $s(t, \alpha)$ $j_k(\alpha)$. (C1 exponent) $\varepsilon_0$. (Conjecture) $C'$.
**(4.2)** (Sums of series) A, $C°$. (C2 exponent) $\varepsilon_1$.
**(4.3)** (Sum of series) $B°$. (Partial fraction expansions) $p(s)$, $p_i(s)$. (C3 exponent) $\varepsilon_2$.
**(4.4)** (C4 exponent) $\tilde{\varepsilon}_1$. (Compound conjecture) $C^\wedge$.
**(4.5)** $v_0$.

**§5** (Cosine sum) $\lambda(y)$. $g_0(y)$. (Sum of exponentials) $e(z)$.

**§6** $j(u)$. (Poisson integral) $\upsilon(z)$. (Sums of exponentials) $e(z, n)$, $\hat{e}(z, n)$. $p_{i,+}(s)$. $\Xi(s)$. $\Theta(\theta, z)$.

**§8**
**(8.2)** **(**Ramanujan tau Dirichlet function) $r(s)$.

**Appendix** (Projection) $P_N(a)$. (Sequence) $\omega$. $S_w(a, J)$. (Abelian group) T. (Metric) $d_N$. $g_k(\theta)$, $G_N(y)$ and $G(y)$. $V_N(y)$ and $D_N(y', y)$.

## Part II

**§1**
$n(s, \alpha)$, $f(s, \alpha)$. $m(z, \beta)$. $Q_k(z)$. $M(x)$. $q(u, \beta)$, $l_k(y, \beta)$.

**§2**
$P(z)$, $g(r, j, q)$, $g(j, q)(r)$, $g(r, j, q, p)$. $q_0(T)$, $j_0(T)$. $\theta_{T,p}(r)$, $\omega_{T,p}(r)$.

**§3**

$n_0(s)$, $f_0(s)$, $h^{<\omega>}$. $N(z, \beta)$, $F(z, \beta)$. $j(u, m)$, $E(v, m)$ $N(u, \beta, 1)$, $F(u, \beta, 1)$. $E(u, \beta, 1)$, $E(u, \beta, 2)$, $J(z)$. $W(z, \beta)$, $B_0(z, \beta)$, $M(z, \beta)$. $I(p, z, u)$, $I(p, z)$.

## §4
$\varphi(B, z)$, $\gamma(\beta, z, *)$, (confluent hypergeometric functions) $M(a, B, z)$. (Kummer's operator) $K(a, B)$. $U(a, a, z)$

## §5
$H(z, \beta)$. $\Omega(r, \delta)$. $A(z, \beta, n)$, $R(z, \beta)$. $G(z, \beta, m)$, $H(z, \beta, 2w)$.

## §6
$n_0(s, \beta)$, $f_0(s, \beta)$, $b(s, \beta)$, $T_0(z, \beta)$, $T_0(z, \beta, 4w)$. $n(s, \beta)$, $f(s, \beta)$, $c(z, \beta)$, $P_{4w}(z, \beta)$. $m(t)$, $d(t_1, t_2)$.

## Part III

**Introduction $\Delta(s)$.** (Complex domain) $S$.

## §1
$B(d)$, $z(s)$, $I(s)$, $G'(s)$, $B_r(d)$, $w(s)$, $T(s)$, $B_r'(d)$, $F_n(z_0, h, z_0)$, $M(h, z_0, \rho)$, $\Delta(s, z)$.

## §2
$B_i(\alpha)$, $k(s)$, $T(s)$, $B_i'(\alpha)$, $Z_i$. $D(t_0, K)$, $S'(k)$.

## Appendix

$E(s)$. $h(z, u, r)$. $\rho(z, u, n)$. (Half-plane) $H_+$. $\rho(z, u)$.

## Part V

$h^{\#}(z)$, $L$, $D(r, \omega)$, $D_1(r, \omega)$, $S(z, \omega)$, $S_{\leq}(z, \omega)$, $\Delta(z, \omega)$, $c(h, k, \omega)$ and $c(h, k)$.

# On the Riemann zeta-function, Part I: Outline

## Introduction

$C^\wedge$ (Part I, §4, (4.4)) is the conjunction of the Riemann hypothesis, RH, (Part I, §2, (2.1)) and hypotheses advanced by the author (Part I, §4, (4.1)-(4.4)). $C^\wedge$ implies the simple zeros conjecture, SZC, (Part I, §2, (2.2), §4, (4.3)). Let $z(n)$ be the nth nonreal ("critical") zero of $\zeta(s)$ having positive imaginary part in order of magnitude thereof. $C^\wedge$ is used to obtain a relation (Part I, §6, Conditional corollary 6.3 (2)) of the pair $z(n)$, $\zeta'(z(n))$ to the preceding such pairs and the $\zeta(\frac{1}{2} + 4w)$ with w a nonnegative integer.

That relation is derived from the two-sided Laplace transform representation of an odd meromorphic function f(s) (Part I, §1) constructed from ζ(½ + s). That representation, for the strip $V_0'$ of s with ½ < Re(s) < 4, (Part I, §3, Main unconditional theorem (1)) is obtained in Part II without using any unproven hypothesis. Assuming $C^\wedge$ , in Part IV the Laplace representation of f(s), on the strip $V_0$ of s with 0 < Re(s) < 4, is derived from the partial fraction expansion of f(s) (Part I, §5, Introduction, Conditional theorems 5.1-5.2) obtained in Part III. The Laplace densities of the unconditional and of the conditional representation of f(s) restricted to $V_0'$ are then shown in Part IV to be equal (Part I, §5, (5.4), Main conditional theorem (1) (i), (ii)). That engenders the above relation for z(n), ζ'(z(n)) in Part V. The functions λ(y) and e(z) (Part I, §5, Introduction), which play a central role therein, are shown (in Part IV, Part V respectively), using $C^\wedge$, to have analytic extensions to the entire complex plane (Part I, §5, (5.4), Conditional theorem 5.3 (1), §6, Conditional corollary 6.3 (1), respectively). It also results from $C^\wedge$ in Part V that a partial fraction expansion, $p_{i,+}(s)$, associated with the roots z(n), is represented, for s with Re(s) < 0, by an integral transform with kernel arising from the above ζ(½ + 4w). See Part I, §6, Conditional corollary 6.5.

In Part IV the hypotheses $C^\wedge$ and C5 (Part I, §4, (4.5)) are used to prove that the Laplace density of f(s) on the strip $V_0$ is positive (Part I, §7, Main conditional theorem (2)). Together that and an unconditional result (Part I, §3, Main unconditional theorem (4)) proven in Part II are used in Part IV to establish that f(s) is a meromorphic characteristic function on the complex plane (Part I, §7, Main conditional theorem (3)). It is shown in Part I, §7, that conversely, RH and SZC hold, if f(s) is an analytic characteristic function on $V_0$. See A. Csizmazia [8-11]. A geometric consequence (Part I, §3, and §7) is obtained in Part VI from a generalization of the Fourier metric geometries of J. von Neumann and I. J. Schoenberg. See A. Csizmazia [12].

In Part I, §8, an analogous context is adumbrated for the recondite case of certain functions built from Dirichlet L-functions, L(s, χ), instead of the zeta-function. The context deriving from the Ramanujan tau Dirichlet function, r(s), is also considered in §8. An introduction to L(s, χ), r(s) is presented in T. M. Apostol [3], [4] respectively.

## §1 Definitions of l(s), a(s), ξ(s), n(s), f(s), b(s), V(x₀, x₁), V[x₀, x₁], V(ε).

Let s be complex. Define:

$$l(s) := \pi^{-s/2} \cdot s\Gamma(s/2) = \pi^{-s/2} \cdot 2\Gamma(1 + s/2), \ a(s) := l(s)(s - 1), \ \xi(s) := (½)a(s)\zeta(s),$$

$\zeta$(s) is the Riemann zeta-function. (See: T. M. Apostol [3], H. M. Edwards [14], A. Ivĭc [18], and E.C. Titchmarsh [39].) Also

n(s) := sin($\pi$s/4)·2$\xi$(½ + s), f(s) :=1/n(s) and b(s) := sin($\pi$s/4)a(½ + s).

Say $x_0 < x_1$. Let V($x_0$, $x_1$) be the open vertical strip of all s with $x_0 <$ Re(s) $< x_1$. Define V[$x_0$, $x_1$] to be the closed strip of all s with $x_0 \le$ Re(s) $\le x_1$. Set V($\epsilon$) := V(0, $\epsilon$) for $\epsilon > 0$.

x, t, y are real variables. k, w, n are integers, except for the function n(s).

**The functional equation of $\zeta$(s). Symmetries of $\zeta$(s), $\xi$(s), n(s), f(s).**
$\xi$(s) is an entire function. $\zeta$(s) satisfies the functional equation

$\xi$(½ - s) = $\xi$(½ + s).

Thus the entire function n(s) is odd: n(-s) = n(s). Hence the meromorphic function f(s) is odd. g(s*) = (g(s))* for g any of $\zeta$, $\xi$, n, f.

Aspects of complex analysis used in Parts I-V are presented in L. Ahlfors [1].

The role of the factor sin($\pi$s/4) in f(s) := 1/(sin($\pi$s/4)2$\xi$(½ + s)) is clarified in Part II, §1.

**Zeros of $\zeta$(s), $\xi$(s), n(s). Poles of f(s).** The zeros of $\xi$(s) are precisely the nonreal zeros of $\zeta$(s). Those zeros give n(z - ½) = 0. The order of the zero is conserved under that translation. Those are the only nonreal zeros of n(s). The real zeros of n(s) arise from . Those zeros are the multiples of four and are simple: n(4w) = 0, n′(4w) ≠ 0 for w = 0, ±1, ±2,... The poles of f(s) correspond to the zeros of n(z) in location and order.

**The nonreal zeros of $\zeta$(s).** These zeros all lie in the "critical" strip V(0, 1) of s with 0 < Re(s) < 1. They are symmetric about each of the real axis and the "critical" line of s with real part ½. G. H. Hardy and J. E. Littlewood [16] proved that an infinite number of the nonreal zeros lie on the critical line. The proportion p(T) of z with $\zeta$(z) = 0 and 0 < Im(z) < T for which also Re(z) = ½ and z is a simple root, $\zeta$′(z) ≠ 0, is ultimately at least 2/5: liminf $_{T\to\infty}$ p(T)≥ 2/5. See: A. Selberg [38], N. Levinson [24], and J. B. Conrey [7].

**§2 Standard conjectures**

**(2.1) The Riemann Hypothesis, RH**

In 1859 B. Riemann [36] formulated the following conjecture.
RH: The real part of each nonreal zero of $\zeta(s)$ is one-half.

The Riemann hypothesis has not been resolved since its formulation despite the determined efforts of generations of leading mathematicians. The intrinsic relation of the zeta function with the primes expressed in its Euler factorization has led to the reduction of the solution of deep questions in multiplicative number theory, in particular concerning the distribution of the primes, to the resolution of RH. The Riemann hypothesis is one of the most celebrated unsolved problems of mathematics. See E. Bompieri [5] and J. B. Conrey [6].

## (2.2) The simple zeros conjecture, SZC.

RH is allied with the unresolved conjecture stated next.

SZC: Each nonreal zero $z$ of $\zeta(s)$ is simple, $\zeta'(z) \neq 0$.

**Definition of $\gamma_n$.** Let $\gamma_1 < \gamma_2 < \ldots < \gamma_n < \gamma_{n+1} \ldots$ enumerate in order of magnitude the distinct imaginary parts $\gamma_n$ of the zeros $z$ of $\zeta(s)$ with $\text{Im}(z) > 0$. When $n \geq 1$, set $\gamma_{-n} = -\gamma_n$.

Computer calculations have established that RH and SZC both hold for the initial ten trillion zeros of $\zeta(s)$ with $\text{Im}(z) > 0$ : If $\zeta(z) = 0$ and $\text{Im}(z) = \gamma_n$ for some $n$ with $1 \leq n \leq 10^{13}$, then $z = \frac{1}{2} + i\gamma_n$ and $\zeta'(z) \neq 0$. See X. Gourdon [15]. See also A. M. Odlyzko [30-34], J. van de Lune, J. J. te Riele and D.T. Winter [26], H. J. J. te Riel and J. van der Lune [35] and H.L. Montgomery [27].

The main conditional theorem of this article assumes RH, SZC and conjectures introduced by the author. That theorem implies RH and SZC. See Part I, §7.

## (2.3) The Lindelof hypothesis, LH.

The Riemann hypothesis implies the Lindelof hypothesis stated next.

LH: If $\varepsilon > 0$, then $|\zeta(\frac{1}{2} + it)| \leq t^\varepsilon$ and large positive $t$.

LH implies that for $\sigma \geq \frac{1}{2}$ and $\varepsilon > 0$: $\zeta(\sigma + it) = O(|t|^\varepsilon)$, for real $t$ with $|t|$ large.

See H. M. Edwards [14] and E.C. Titchmarsh [39].

## (2.4) The Lehmer phenomenon

The "Lehmer phenomenon" is that certain "Lehmer pairs" of successive zeros of ζ(s) on the critical line of s with Re(s) = ½ are exceptionally close.

**Definition of $\delta_k$.** Set $\delta_k := \min\{\gamma_k - \gamma_{k-1}, \gamma_{k+1} - \gamma_k\}$.

The conjectured Gaussian unitary ensemble (GUE) model (See J. P. Keating and N. C. Smith [19-21]) indicates and calculations confirm that there are Lehmer pairs ½ + i$\gamma_k$ , ½ + i$\gamma_{k+1}$ for which the values of $\delta_k$ , |ζ′(½ + i $\gamma_k$)| and |ζ(½ + im)| , with m the midpoint ($\gamma_k$ + $\gamma_{k+1}$ )/2, can be as low in order as $\gamma_k^{-\frac{1}{3}}$, $\gamma_k^{-\frac{1}{3}}$ and $\gamma_k^{-2/3}$ respectively. See D. H. Lehmer [23]; H. M. Edwards [14], Lehmer's Phenomenon, §8.3 pp. 175-179; G. Csordas, W. Smith and R. S. Varga [13]; A. M. Odlyzko [29]; X. Gourdon [15]; and Eric W. Weisstein [41].

This appears as the worst case scenario which renders the proof of the main conditional theorem of this article more delicate.

## §3 Context and statement of the main unconditional theorem

A result is termed "unconditional" when it is (or can be) established without relying on unproven conjectures.

**Definitions of $V_u$, $V_u′$.** Say u is a multiple of four, u = 4w. Let $V_u := V(u, u+4)$. If u ≠ 0, -4, set $V_u′ = V_u$. Let $V_0′ := V(½, 4)$. Take $V_{-4}′ = -V_0′$.

f(-s) = f(s). So assume w ≥ 0. f(s) is analytic on $V_{4w}′$

**Integrability of |f(s)| on vertical lines outside of the critical strip.**

**Theorem 3.1** *|f(s)| is both integrable and square integrable on each vertical line x + i$\boldsymbol{R}$ with x ≥ ½ and x ≠ 4w for w ≥ 1.*

f(s*) = (f(s))*. So Theorem 1 holds if $\int_{t>0}(dt)|f(x + it)|^p < \infty$ for p = 1, 2.

Theorem 3.1 is proven in Part II, §3.

We determine the representation for s in $V_{4w}′$ of $(-1)^w f(s)$ as a two-sided Laplace transform: $(-1)^w f(s) = \int_{\boldsymbol{R}} d(y)e^{sy}P_{4w}(\pi e^{-2y})$. The equivalent Fourier transform representation is obtained by the rotation s = iz. The equivalent Mellin transform representation is $(-1)^w f(s) = \int_{v>0} d(v)v^{s-1}P_{4w}(\pi \cdot v^{-2})$.

f(s*) = (f(s))* implies that each measure density $P_{4n}(v)$ is real for v > 0. f(-s) = f(s) implies that for w ≤ -1: $P_{4w}(\pi v) = P_{-4(w+1)}(\pi/v)$ for v > 0.

**Definition of the Pochhammer symbol $(z)_n$.** $(z)_n := \Pi_{0 \le k \le n-1}(z + k)$, with z complex and the integer $n \ge 1$. Take $(z)_0 = 1$.

**Definitions of $\tilde{c}(4k)$, $c(z)$.** Let k be an integer $\ge 0$.

$$\tilde{c}(4k) := 1/(\pi^{\frac{3}{4}}\Gamma(5/4 + 2k)(2k - \tfrac{1}{4})\zeta(\tfrac{1}{2} + 4k)).$$

$\tilde{c}(4k) > 0$. Define $c(z) := 1/n'(z)$, for z with $n'(z) \ne 0$. $c(4k) = \tilde{c}(4k)(-(\pi^2))^k$. In particular $c(0) = 2^4/(\pi^{\frac{3}{4}}\Gamma(\tfrac{1}{4})\cdot(-\zeta(\tfrac{1}{2})))$.

The $\tilde{c}(4k)$ converge ultrarapidly to 0 as $k \to \infty$. This is a result of the following observations (i), (ii).
(i) Say $p > 1$. Comparing $\zeta(p)-1$ to $\int_{x > 2} (dx)x^{-p}$ gives $1 + 1/(2^{p-1}\cdot(p-1)) < \zeta(p) < 1 + 2^{-p}(2/(p-1) + 1)$
(ii) The Stirling approximation to $\Gamma(x)$, for $x > 0$, leads to

$$1 / \tilde{c}(4k) \sim (2\pi)^{\frac{1}{2}} (2k)^{7/4}(2k/e)^{2k}(1 + O(1/k)).$$

**Definition of $P_0(z)$.** Set $P_0(z) := (-1)\sum_{k \ge 1} \tilde{c}(4k)(-(z^2))^k$.

Then $P_0(z)$ is an entire function of z. $P_0(z)$ is even: $P_0(-z) = P_0(z)$. Also $P_0(z^*) = (P_0(z))^*$.

**Theorem 3.2** (i) *$P_0(z) \sim O(z^2)$ for z near 0.*
   (ii) *Say $\varepsilon > 0$. Then for $v \ge 0$, $|P_0(v)| \sim O(v^{\frac{1}{4}+\varepsilon})$, as $v \to \infty$.*

Theorem 3.2 (ii) is established in Part II, §6.

Part (i) of the previous theorem has the following consequence. Let y be real.
If $x < 4$, then $e^{xy}|P_0(\pi e^{-2y})| \sim O(e^{-(4-x)y})$ as $y \to \infty$.
Part (ii) implies the following.
If $x > \tfrac{1}{2}$, then for any $\varepsilon$ with $0 < \varepsilon < x - \tfrac{1}{2}$, $e^{xy}|P_0(\pi e^{-2y})| \sim O(e^{\varepsilon y})$ as $y \to -\infty$.
Thus ('): $\tfrac{1}{2} < \text{Re}(s) < 4$ implies $|e^{sy}P_0(\pi e^{-2y})|$ vanishes with exponential rapidity as $|y| \to \infty$.

The nonzero coefficients in the Taylor expansion of $P_0(z)$ alternate in sign. That is an obstacle to using the expansion to analyze the behavior of $P_0(v)$ as $v \to \infty$ on the real line. That applies to $P_0(\pi \cdot e^{-2y})$, when $y \to -\infty$ thereon.

$P_0(v)$ is bounded on the real axis, provided $C^\wedge$ holds. That is stated in §5, (5.4), as Main conditional theorem (1) (i'), after presenting the compound conjecture $C^\wedge$ in §4, (4.4). It is proven in Part IV, §2.

**Theorem 3.3** *$P_0(v)$ is strictly increasing with $v$ for $0 \le v \le \pi$. $P_0(0) = 0$.*

The previous theorem is proven in Part II, §6,

$P_0(v) > 0$, if $0 < v \le \pi$ (or $-\pi \le v < 0$). Later we prove the conditional result that $P_0(v) > 0$ for $v > \pi$. (See C5 of §5, (4.5), Main conditional theorem (2) of §7 herein and of § 3 of Part IV.) Then $P_0(v) > 0$ for all nonzero $v$.

The next lemma is used to establish the convergence and analyticity properties of certain essential partial fraction representations.

**Definition of the open disk B(z, r).** Say $r \ge 0$. Take $B(z, r) := \{s: |s - z| \le r\}$.

The previous connotation of $B(z, \beta)$ and that for the complex function $B(z, r)$ of Part II, §3, are to be identified by the context.

Let S be a subset of the complex plane C. Fix the positive integer p. Set $\omega_p :=$ $\exp(i2\pi/p)$. Let $S_p := U_{0 \le k \le p-1} \omega_p^k \cdot S$. Say $B(0, r) \cap S$ is finite for each r. Let $t(s) \ge 0$, for each s in S. $\sum_{s \varepsilon S} t(s) := \lim_{r \to \infty} \sum_{s \varepsilon B(0, r)} t(s)$. Now let $t(s)$ be complex. Assume $\sum_{s \varepsilon S} |t(s)|$ is finite. Then $\sum_{s \varepsilon S} t(s) := \lim_{r \to \infty} \sum_{s \varepsilon B(0, r) \cap S} t(s)$.

Assume Z is an infinite subset of C, with $B(0, r) \cap Z$ finite for each r. Say for any distinct z, z′ in Z, $z \ne \omega_p^k \cdot z′$, for $1 \le k \le p - 1$. Let $c_z$ be a nonzero complex number, for z in Z. Take $\hat{G}(s) := \sum_{z \varepsilon Z} |c_z/(s^p - z^p)|$, for s in C - $Z_p$. Say $r \ge 0$ and $d > 0$. Let $K(r, d)$ be the set of all s with $|s| \le r$ and $|s - z| \ge d$, for each z in $B(0, r + d) \cap Z_p$. If s is in $K(r, d)$, then $|s - z| \ge d$, for each z in $Z_p$. $K(r, d)$ is a compact subset of C - $Z_p$. We will use the cases p = 1, 2 of the next lemma.

**Lemma 3.1** *Say $\hat{A} := \sum_{z \varepsilon Z, z \ne 0} |z|^p |c_z|$ is finite.*
*(1) Let K be a compact subset of C - $Z_p$. As $\rho \to \infty$, $\sum_{z \varepsilon B(0, \rho) \cap Z} |c_z/(s^p - z^p)|$ converges to a finite value $\hat{G}(s)$, uniformly in s on K. So $G(s) := \sum_{z \varepsilon Z} c_z/(s^p - z^p)$ behaves likewise.*
*(2) $\hat{G}(s)$ is continuous on C - $Z_p$. $G(s)$ is analytic on C - $Z_p$. $G(s)$ has a simple pole at each z in $Z_p$.*

**Proof of (1)**. Say K is non empty. K is contained in $K(r, d)$, with r, d as follows. Set $r = \max\{|s|: s \varepsilon K\}$. Let $d′ = \min\{|z|: z \varepsilon Z$ and $|z| > r\}$ - r. If $B(0, r) \cap Z$ is empty, set $d = d′$. Say $B(0, r) \cap Z$ is non empty. Let $d″$ be the minimim distance between K and $B(0, r) \cap Z_p$. Take $d = \min\{d′, d″\}$.

Claim
(1) *Restrict s to K(r,d). Say z is in Z. $|c_z/(s^p - z^p)| \le t(z)$, with $t(z)$ as follows. If z*

*is in $B(0, r + d) \cap Z$, then $t(z) := d^p|c_z|$. If $z$ is in $Z$ and $|z| > r + d$, then $t(z) := j^{-1} \cdot |z|^{-p}|c_z|$, with $j = 1 - (r/(r + d))^p$.*

(2) $\sum_{z \in Z} t(z) = d^{-p} \cdot \sum_{z \in Z, \, |z| \leq r + d} |c_z| + j^{-1}(\hat{A} - \sum_{z \in Z, \, z \neq 0, \, |z| \leq r + d} |z|^{-p}|c_z|) < \infty.$

**Proof of Claim (1).** Say $z$ is in $B(0, r + d) \cap Z$. $|s^p - z^p| = \Pi_{0 \leq k \leq p-1} |s - \omega_p^k \cdot z| \geq d^p$. Suppose $z$ is in $Z$ and $|z| > r + d$. $z$ is nonzero. $|z|^{-p}|s^p - z^p| = |1 - (s/z)^p 1| \leq 1 - (|s| / |z|)^p$. Also $|s| \leq r$. Thus $|z|^{-p}|s^p - z^p| \geq j$.

**Definition of N.** Let **N** be the set of positive integers.

**Definition and properties of F(s).**

**Corollary 3.1** *$\sum_{w \geq 1} |c(4w)/(s - 4w)|$ converges. $F(s) := \sum_{w \geq 1} c(4w)/(s - 4w)$ is analytic on $C - 4N$. Also $F(s)$ has a simple pole at each $4w$ with $w \geq 1$.*

**Definitions** $p_r(s)$, **Z**, $Z_i$, $Z°$. Set $p_r(s) := (c(0)/s) + F(s) - F(-s)$. $p_r(s)$ is analytic except for a simple pole at each $4w$, with $w$ an integer. Let $Z$ be the set of integers. Let $Z_i$ be the set of nonreal zeros of $\zeta(\frac{1}{2} + s)$. Set $Z° = (4Z)UZ_i$.

**Main unconditional theorem**

**(1)** *On the strip $V_0'$ of $s$ with $\frac{1}{2} < Re(s) < 4$: $f(s)$ is analytic and*

$$f(s) := 1/(\sin(\pi s/4) \cdot 2\xi(\frac{1}{2} + s)) = \int_R d(y)e^{sy}P_0(\pi e^{-2y}), \text{ with}$$

$$P_0(z) := (4/(\pi^{\frac{1}{4}}\Gamma(\frac{1}{4})) \cdot (-1)\sum_{k \geq 1}(-(z^2))^k/((5/4)_{2k}(2k-\frac{1}{4})\zeta(\frac{1}{2} + 4k)) \text{ entire.}$$

**(2)** *$f(s) - p_r(s)$ has an analytic extension from $C - Z°$ to $C - Z_i$. The open set $C - Z_i$ includes the $s$ with $|Re(s)| \geq \frac{1}{2}$.*

**(3)** *$F(s) = \int_{y > 0} d(y)e^{sy}P_0(\pi e^{-2y})$, for $s$ with $Re(s) < 4$.*

**Proof of (3).** Let $x := Re(s)$. $1/s = \int_{y > 0} d(y)e^{sy}(-1)$, when $x < 0$. Say $x < 4$. Then $F(s) = \sum_{w \geq 1}\int_{y > 0} d(y)t(w, s, y)$, with $t(w, s, y) := e^{sy}(-1)c(4w)e^{-4wy}$. Now $|t(w, s, y)| \leq |c(4w)| \cdot e^{(x-4)y}$. So $\int_{y > 0}\sum_{w \geq 1}|t(w, s, y)|$ is finite. Therefore one has $\sum_{w \geq 1}\int_{y > 0} t(w, s, y) = \int_{y > 0}\sum_{w \geq 1} t(w, s, y)$.

Together (1) and (3) imply that $f(s) = F(s) + \int_{y < 0} d(y)e^{sy}P_0(\pi e^{-2y})$, on $V_0'$. Main unconditional theorem (3) yields the following corollary.

**Corollary 3.2** *$p_r(s) = \int_{y < 0} d(y)e^{sy}(c(0) - P_0(\pi e^{2y})) + \int_{y > 0} d(y)e^{sy}P_0(\pi e^{-2y})$, provided $0 < Re(s) < 4$.*

The inverse two-sided Laplace transform, $g(y)$ of $f(s)$ on $V_0'$ is determined in the Main unconditional theorem by the explicit formula $g(y) = P_0(\pi e^{-2y})$. This formula has the following noteworthy aspects.

(1) It holds even when y < 0.

(2) The Taylor series for $P_0(z)$ has coefficients which apart from elementary expressions involve only the values $\zeta(\frac{1}{2} + 4k)$ with $k \geq 1$. See A key observation, Part I, §6. Despite the relations $f(s) := (1/b(s))\cdot(1/\zeta(\frac{1}{2} + s))$ and $1/\zeta(u) = \sum_{n \geq 1} \mu(n) / n^u$ provided $Re(u) > 1$, the Moebius function $\mu(n)$ does not directly appear in the Taylor series for $P_0(z)$.

(3) The Taylor series shows that $P_0(z)$ extends from the ray z > 0 to the complex plane C as an even entire function.

**Definition of g(z).**

(4) $g(z) := P_0(\pi e^{-2z})$ extends g(y) to an entire function with period $i\pi/2$: $g(z + i\pi/2) = g(z)$.

Next we consider f(s) on $V_{4w}$ with w > 0. f(s) is analytic on $V_{4w}$.

**Definition of $P_{4w}(z)$.** Say w = 1,2,3,... Let

$$P_{4w}(z) := (-1)^w(P_0(z) + \sum_{1 \leq k \leq w} \tilde{c}(4k)(-(z^2))^k) = (-1)^{w+1}\sum_{k \geq w+1}\tilde{c}(4k)(-(z^2))^k.$$

$P_{4w}(z)$ is an even entire function of z. $P_{4w}(0) = 0$. $P_{4w}(z) \sim O(z^{2(w+1)})$ for z near 0. Say $v \geq 0$ and $v \to \infty$. Say $\varepsilon > 0$. Then $P_4(\pi v) = |c(4)| v^2 + O(v^{\frac{1}{4}+\varepsilon})$. Say $w \geq 2$. Then $P_{4w}(\pi v) = |c(4w)| v^{2w} + O(v^{2(w-1)})$. Hence for $w \geq 1$, $P_{4w}(v) \to \infty$ as $v \to \infty$ on the real line.

Say y is real. If x > 4w, then as $y \to -\infty$, $e^{xy}|P_{4w}(\pi e^{-2y})| \sim O(e^{(x-4w)y})$. If x < 4(w+1), then as $y \to \infty$, $e^{xy}|P_{4w}(\pi e^{-2y})| \sim O(e^{-(4(w+1)-x)y})$. Thus $4w < Re(s) < 4(w+1)$ implies that $|e^{sy}P_{4w}(\pi e^{-2y})|$ decays to zero with exponential rapidity as $|y| \to \infty$.

**Lemma 3.2** *Fix $w \geq 1$.*

*(i) $P_{4w}(v)$ is strictly monotone increasing in v on the interval $0 \leq v \leq \pi$. $P_{4w}(0) = 0$.*

*(ii) $P_{4w}(v) > 0$ for v > 0.*

Lemma 3.2 and the following theorem are proven in Part II, §6.

**Main unconditional theorem**

**(4) *Say $w \geq 1$.***

**(i) *On $V_{4w}$, f(s) is analytic and***

$$(-1)^w f(s) = \int_R d(y)e^{sy}P_{4w}(\pi e^{-2y}),$$

***with $P_{4w}(z)$ entire.***

**(ii) *$P_{4w}(\pi e^{-2y}) > 0$ for all real y.***

**The alteration of the density under a transition of domain of a Laplace representation.**

The principle presented next predicts the change in the Laplace density of f(s) induced by the transition of the domain of s in the unconditional case from $V_0'$ to $V_4$ or from $V_{4w}$ to $V_{4(w+1)}$, with $w \geq 1$; and in the conditional case, using the Main conditional theorem of Part I, §5, (5.4), from $V_0$ to $V_4$.

Say N(s) is analytic on the open subset D of the complex plane. Assume N(r) = 0 and N'(r) ≠ 0. Let $F_1(s) = F_1(s, N, r)$ be N(s)/(s-r) for s ≠ r and $F_1(r) := N'(r)$. Let $F_2(s) = F_2(s, N, r)$ be $(F_1(s) - N'(r))/(s-r)$ for s ≠ r and $F_2(r) := N''(r)/2$. Then $F_1(s)$, $F_2(s)$ are analytic on D. On D the zeros of $F_1(s)$ are precisely the zeros z of N(s) other than r and the order of the zero z is conserved. $F_2(s)$ does not vanish at those z. Also 1/N(s) = 1/(N'(r)(s - r)) - (1/N'(r))($F_2(s)/F_1(s)$).

Suppose that D, N, r are as above with D = $V(x_0, x_1)$, $x_0 < r < x_1$ and r the only zero of N(s) on $V(x_0, x_1)$. Assume (1), (2) below.
(1) $\int$d(y)$e^{xy}$|j(y)| is finite, respectively with y > 0, y < 0 and respectively $x_1 - x$ , x - $x_0$ positive and arbitrarily small,
(2) If s is on $V(x_0, x_1)$, then -(1/N'(r))($F_2(s)/F_1(s)$) = $\int_R$d(y)$e^{sy}$j(y).
Define $m_0(y) := j(y) - (e^{-ry}/N'(r))$ for y > 0 and $m_0(y) := j(y)$ for y < 0.
Set $m_1(y) := j(y)$ for y > 0 and $m_1(y) := j(y) + (e^{-ry}/N'(r))$ for y < 0. Then (1′), (2′) below hold.
(1′) If s is on $V(x_0, r)$, respectively $V(r, x_1)$ then $\int_R$ d(y)|$e^{sy}m_k(y)$| is finite and 1/N(s) = $\int_R$ d(y)$e^{sy}m_k(y)$ with k respectively 0, 1.
(2′) $m_1(y) = m_0(y) + (e^{-ry}/N'(r))$.
Conversely, assume (1′) and (2′). Let j(y) be defined as $m_0(y)$ for y > 0 and $m_1(y)$ for y < 0. Then (2) holds with $\int_R$ d(y)|$e^{sy}$j(y)| finite.

**A geometric consequence of the Main unconditional theorem (4).**

**Definitions of $m_x(t)$, $d_x(t_1, t_2)$.** Say x, t are real and x is not a multiple of four. Set $m_x(t) := |1 - (n(x)/n(x + it))|^{1/2}$ , if n(x + it) is nonzero. (Take $m_x(t) = \infty$, if n(x + it) = 0.) Let $d_x(t_1, t_2) := m_x(t_1 - t_2)$, for real $t_1$, $t_2$.

The next assertion is proven in Part VI. That proof uses the Main unconditional theorem (4).

*Fix x with: either x < -4, or x > 4; and with x not a multiple of four.*
*$m_x(t)$ is a (finite-valued) metric norm in t on the real line.*
*$d_x(t_1, t_2)$ is a translation invariant metric in $t_1$, $t_2$ on the real line.*

## §4 Conjectures introduced by the author

Applying Theorem 3.1 with s on $V_0'$, Theorem 3.2 (i) and the Main unconditional theorem of §3 yields the following corollary.

**The order of $P_0(v)$ and the partial validity of RH.**

**Conditional corollary 4.1** *Assume that:*
(*) $0 \leq j < \frac{1}{4}$ *and for arbitrarily small $\varepsilon > 0$, $|P_0(v)| \sim O(v^{j+\varepsilon})$, for $v \geq 0$ and $v$ large.*
*Then:*
(1) $f(s) = \int d(y)e^{sy}P_0(\pi e^{-2y})$ *on all of the strip $V(2j, 4)$.*
(2) $\zeta(s)$ *does not have any zeros $z$ with $\frac{1}{2} + 2j < Re(z) < 4$. If $j = 0$, then the Riemann hypothesis is valid.*

**Proof** The assumption (*) implies $\int_{y<0} d(y)e^{xy}|P_0(\pi e^{-2y})|$ is finite, when $x > 2j$. Therefore $\int_{y<0} d(y)e^{sy}P_0(\pi e^{-2y})$ converges to an analytic function on the half-plane of s with $Re(s) > 2j$. Then $E(s) := \int d(y)e^{sy}P_0(\pi e^{-2y})$ is analytic, provided $2j < Re(s) < 4$. Now $E(s) = f(s)$, when $\frac{1}{2} < Re(s) < 4$. Therefore $E(s) = f(s)$ also when $2j < Re(s) \leq \frac{1}{2}$. So (1) and (2) hold. $\square$

Assume RH. Then $f(s)$ is analytic on $V(\frac{1}{2})$ and hence also on $V_0$, $V_{-4}$.

We assume RH and SZC in establishing a representation $f(s) = \int_R d(y)e^{sy}g_0(y)$ valid on the entire strip $V_0$. An explicit formula for $g_0(y)$ is determined, which for $y < 0$ involves the imaginary parts, $\gamma_k$ for $k \geq 1$, of the critical zeros of $\zeta(s)$ in the upper half-plane. It also involves the values of $\zeta'(s)$ at those zeros.

The derivation of this formula is accomplished by assuming, in addition to RH and SZC, the conjectures made by the author and presented next.

## (4.1) Conjecture 1

RH and the first conjecture, C1, will serve to establish that $|f(s)|$ is integrable (and square integrable) on any vertical line in $V(\frac{1}{2})$.

Say $x > 0$ and $x \neq 4w$. The integrability of $|f(x + it)|$ over nonnegative t depends on its asymptotic behaviour as $t \to \infty$. We shall reduce that question to a study of $1/|\zeta(\frac{1}{2} + x + it)|$ as follows.

**The asymptotic behavior of $|b(s)|$ on a vertical strip of finite width**

$f(s) := (1/b(s))\cdot(1/\zeta(\frac{1}{2} + s))$. We obtain the asymptotic behavior of $|b(s)|$ on the vertical strip of $s = x + it$ with $x_0 < x < x_1$, $t$ real, using the Stirling formula for $\Gamma(z)$ thereon for $|t| \geq T > 0$. (See G. Andrews, R. Askey, R. Roy [2].)

$$|\Gamma(x + it)| \sim (2\pi)^{\frac{1}{2}} \cdot |t|^{x - \frac{1}{2}} \cdot e^{-\frac{1}{2}\pi|t|}(1 + \varepsilon(x, t)/|t|), \text{ with } |\varepsilon(x, t)| < K(x_0, x_1) < \infty.$$

$$|b(x + it)| \sim K_1(x)\cdot|t|^{7/4 + x/2}(1 + \varepsilon(x, t)/|t|),$$

$$K_1(x) := (\pi/2)^{\frac{1}{4}} \cdot (2\pi)^{-x/2}, \; 0 < K_1(x_1) < K_1(x) < K_1(x_0) \text{ and } |\varepsilon(x, t)| < K_2(x_0, x_1).$$

The integrability of $|f(x + it)|$ over $t \geq 0$ when $0 < x < \frac{1}{2}$ is later deduced from RH and the Conjecture 1′ formulated below after some preliminary definitions.

**Definitions of $\delta_k'$, $I_k(\alpha)$, $S_k(\alpha)$, $T_*(\alpha)$, $t(1)$, $x(t, \alpha)$, $s(t, \alpha)$, $j_k(\alpha)$.**
Assume RH. Set $\delta_k' := \min\{1/\log(|\gamma_k|), \gamma_k - \gamma_{k-1}, \gamma_{k+1} - \gamma_k\}$. Say $0 < \alpha < \frac{1}{2}$. Define $I_k(\alpha)$ to be the open interval $\{t\colon \gamma_k + \alpha\delta_k' < t < \gamma_{k+1} - \alpha\delta_{k+1}'\}$, $S_k(\alpha)$ to be the semicircle $\{s\colon |s - i\gamma_k| = \alpha\delta_k', \text{Re}(s) \geq 0\}$ and $T_*(\alpha) := U_{k \geq 1}\, (I_k(\alpha) \cup S_k(\alpha))$. Set $t(1) := \gamma_1 - \alpha\delta_1'$. Say $t > t(1)$. Let $x(t, \alpha)$ be the unique real $x$ with $x + it$ in $T_*(\alpha)$. Set $s(t, \alpha) := x(t, \alpha) + it$. Let $j_k(\alpha) := \min\{|\zeta(\frac{1}{2} + it)/(t - \gamma_k)| : t \text{ is real and } |t - \gamma_k| \leq \alpha\delta_k'\}$, with $\zeta(\frac{1}{2} + it)/(t - \gamma_k)$ at $t = \gamma_k$ interpreted as being $\zeta'(\frac{1}{2} + i\gamma_k)$. One then has $|\zeta'(\frac{1}{2} + i\gamma_k)| \geq j_k(\alpha)$.

**C1′ = Conjecture 1′** There exists an $\alpha$ with $0 < \alpha < \frac{1}{2}$ such that:

$$\int_{t > t(1)} (dt)(t^{7/4} |\zeta(\frac{1}{2} + s(t, \alpha))|)^{-1} < \infty.$$

We assume the following stronger conjecture.

**C1 = Conjecture 1** Assume $0 < \alpha < \frac{1}{2}$, and (i), (ii) as follows hold.
(i) There exist $\varepsilon$, $\lambda$ with $\varepsilon > 0$, $\lambda > 0$ such that $|\zeta(\frac{1}{2} + s(t, \alpha))| > \lambda \cdot t^\varepsilon$, for $t > t(1)$.
**Definition of $\varepsilon_0$.** Let $\varepsilon_0$ be the infimum of such $\varepsilon$.
(ii) $\varepsilon_0 < \frac{3}{4}$.

GUE suggests and calculations corroborate that for some Lehmer pairs $\frac{1}{2} + i\gamma_k$, $\frac{1}{2} + i\gamma_{k+1}$, the order of $|\zeta(\frac{1}{2} + i(\gamma_k + \gamma_{k+1})/2)|$ can descend as low as $\gamma_k^{-2/3}$.

**Definition of C′.** Let C′ be the assumption that each of RH, C1 (i) and $\varepsilon_0 < 7/4$ holds.

C′ will be used to deduce that $f(x + it)$ approaches 0 as $t \to \infty$, when $x(t, \alpha) \leq x < \frac{1}{2}$. See the Conditional lemma 2.4 of Part III, §2.

**(4.2) Conjecture 2**

Recall that c(z) := 1/n′(z), for z with n′(z) ≠ 0. n′(z) and c(z) are even. Hence n′ and c are real-valued on the imaginary axis. If RH and SZC hold, then $c(i\gamma_k)$ = 1/(b($i\gamma_k$)ζ′(½ + $i\gamma_k$)).

The following series enter in the determination of the Laplace density $g_0(y)$ on the negative ray y ≤ 0:

**Definitions of A, C°.** Set A := $2\sum_{k \geq 1}$ |c($i\gamma_k$)|, C° := $\sum_{k \geq 1}$ |c($i\gamma_k$)|/($\gamma_k^2$) and $\sum_{k \geq 1}$ |c($i\gamma_k$)cos($\gamma_k$y)|. (See §5.)

Those series converge, if that for A does. A < ∞ follows from the next conjecture together with RH as established with Conditional lemma 4.1.

## C2 = Conjecture 2

(i) There exists a real ε (with ε ≥ 0) such that: for any σ > 0, there is a K > 0 with |ζ′(½ + $i\gamma_k$)| > K$\gamma_k^{-(\varepsilon + \sigma)}$, for all k ≥ 1.

**Definition of $\varepsilon_1$.** Let $\varepsilon_1$ be the least ε as in (i).

(ii) $\varepsilon_1$ < ¾.

The conjectured GUE model suggests that $\varepsilon_1$ is of order ⅓.

RH and part (i) of C2 together imply SZC. Together RH and |c($i\gamma_k$)| < ∞ for all k ≥ 1, imply SZC.

It is known that $\gamma_k$ ~ 2πk/log(k) as k → ∞. Say 0 < q < 1. Then for all large k, one has $\gamma_k$ > $k^q$.

**Conditional lemma 4.1** *Assume RH and C2. Then A < ∞.*

**Proof** 1/ |c($i\gamma_k$)| = |n′($i\gamma_k$)|, with n′($i\gamma_k$) = b($\gamma_k$)ζ′(½ + $i\gamma_k$). Stirling's formula for Γ(z) shows that |b(it)| ~ (π/2)^¼ ·|t|^{7/4} (1 + O(1/|t|)), with β > 0 , t real and |t| large. C2 then implies that for any ε > 0, there is a K > 0 such that |c($i\gamma_k$)| ≤ K$\gamma_k^{-(1 + \theta - \varepsilon)}$, for all k ≥ 1, where θ = ¾ - $\varepsilon_1$ > 0. Take q with 1/(1+θ) < q < 1. Take ε with 0 < ε < θ and j > 1, with j := (1 + θ - ε)q. Then |c($i\gamma_k$)| ~ O($k^{-j}$) for large k. Thus A is finite.

**Note on the finiteness of A:** *A priori* it is possible that A is finite even if $\varepsilon_1$ = 7/4. Suppose that for m ≥ 1, |ζ′(½ + $i\gamma_{k(m)}$)| ≥ K·$\gamma_{k(m)}^{-p(m)}$, with p(m) > 0. Set q(m) := 7/4 – p(m). Assume q(m) > 0 and $\lim_{m \to \infty}$ q(m) = 0. Say 0 < ρ < 1 and k is large. $\gamma_k$ > $k^\rho$. |c($i\gamma_{k(m)}$)| ≤ $K_1$·$\gamma_{k(m)}^{-q(m)}$ < $e^{-r(m)}$, with r(m) := ρq(m)log(k(m)). $\sum_{m \geq 1}$ $e^{-r(m)}$ converges, if r(m) → ∞ with sufficient rapidity as m → ∞.

## (4.3) Conjecture 3

**Definition of B°.** Let $B° := \sum_{k \geq 1} |c(i\gamma_k)|/\delta_k'$.

The next conjecture C3 will be employed to prove that $B° := \sum_{k \geq 1} |c(i\gamma_k)|/\delta_k'$ is finite. (See the next Conditional lemma 4.2.) $B° < \infty$ implies A is finite, since $0 < \delta_k' \leq 1/\log(\gamma_k) \leq 1/\log(2)$, for $k \geq 1$. $B° < \infty$ is used in the first step in obtaining the density $g_0(y)$. That step assumes $C^\wedge$ of §4, (4.4). It consists of proving that $f(s) := 1/n(s)$ is represented for all s, other than the zeros of $n(s)$, by its formal Mittag–Leffler partial fraction expansion $p(s)$ defined below. (See §5, Introduction, Conditional theorem 5.1.) That proof is achieved in Part III, §2, Conditional theorem 2.2.

**Definitions of p(s), $p_i(s)$.** Assume RH and that A is finite. Say $n(s)$ is nonzero.

$$p(s) := \sum_{z: n(z) = 0} (1/n'(z))(1/(s - z)).$$

$$p(s) = p_i(s) + p_r(s), \text{ with}$$

$$p_i(s) := \sum_{k \geq 1} c(i\gamma_k)(1/(s - i\gamma_k) + 1/(s + i\gamma_k)) = 2s\sum_{k \geq 1} c(i\gamma_k)(1/(s^2 + \gamma_k^2)).$$

**Conditional claim 4.1** *Say $C° < \infty$. Then $\sum_{k \geq 1} c(i\gamma_k)(1/(s - i\gamma_k) + 1/(s + i\gamma_k))$ is absolutely convergent on $C - Z°$. Also $p_i(s)$ is analytic except for simple poles at $\pm i\gamma_k$ for $k \geq 1$. Assume RH and that A is finite. Then $p(s)$ is analytic except for simple poles $\pm i\gamma_k$ for $k \geq 1$, and at 4w, for w an integer.*

**Proof** Apply Lemma 3.1 to the series for $p_i(s)$.

**Conditional claim 4.2** *Assume RH and $C° < \infty$. Then $\Delta(s) := f(s) - p(s)$ has an analytic extension from $C - Z°$ to C.*

**C3 = Conjecture 3**
(i') There exists a real $\varepsilon$ (with $\varepsilon \geq 0$) such that: for any $\sigma > 0$, there is a $K > 0$ with $|\zeta'(\frac{1}{2} + i\gamma_k)| > K\gamma_k^{-(\varepsilon + \sigma)}$, for all $k \geq 1$.
**Definition of $\varepsilon_1$.** Let $\varepsilon_1$ be the least $\varepsilon$ as in (i).
(i) There exists an $\varepsilon \geq 0$ such that: for any $\sigma > 0$, there is a $K > 0$ with $\delta_k > K\gamma_k^{-(\varepsilon + \sigma)}$, for all $k \geq 1$.
**Definition of $\varepsilon_2$.** Let $\varepsilon_2$ be the least $\varepsilon$ as in (i).
(ii) $\varepsilon_1 + \varepsilon_2 < \frac{3}{4}$.

GUE and calculations indicate that $\varepsilon_1, \varepsilon_2$ are of order $\frac{1}{3}$ and so $\varepsilon_1 + \varepsilon_2 \sim \frac{2}{3}$.

Note that C3 (i') is C2 (i). Also C3 (ii) implies C2 (ii). Thus C3 entails C2. Utilize Conditional lemma 4.1 to obtain the following. If RH and C3 hold, then A is finite.

**Conditional lemma 4.2** *C3 implies B° is finite.*

**Proof** C3 implies C2. In the proof of the Conditional lemma 4.1 of (4.2) replace $\varepsilon_1$ with $\varepsilon_1 + \varepsilon_2$.

## (4.4) Conjecture 4

The next conjecture is also used to establish the partial fraction decomposition $f(s) = p(s)$.

### C4 = Conjecture 4

There exists an $\alpha$ as in C1 for which (i), (ii) as follows also hold.
(i) There is an $\varepsilon \geq 0$ such that: for any $\sigma > 0$, there is a $K > 0$ for which for all $k \geq 1$: $j_k(\alpha) > K\gamma_k^{-(\varepsilon + \sigma)}$.
**Definition of $\tilde{\varepsilon}_1$.** Let $\tilde{\varepsilon}_1$ be the least $\varepsilon$ as in part (i). Assume C3 (i).
(ii) $\tilde{\varepsilon}_1 + \varepsilon_2 \leq 1$.

C4 (i) implies C2 (i). $\tilde{\varepsilon}_1 \geq \varepsilon_1$. C3 (ii) and C4 (ii) together amount to $\max\{\varepsilon_1, \tilde{\varepsilon}_1 - \frac{1}{4}\} + \varepsilon_2 < \frac{3}{4}$.

We will need C3 (ii) and $\frac{3}{4} - (\varepsilon_1 + \varepsilon_2) + 1 - (\tilde{\varepsilon}_1 + \varepsilon_2) > 0$, rather than the stronger C4 (ii). See the proof of the Conditional theorem 2.1 of Part III, §2.

**Definition of C^.** Let C^ be the compound conjecture that RH, C1, C3 and C4 all hold.

Note that C^ implies A is finite. (See **(4.3)**.)

## (4.5) Conjecture 5

**Definition of an analytic characteristic function on a vertical strip.**
Say $-\infty \leq w_0 < w_1 \leq \infty$. Let V be the vertical strip of $s = x + it$ with $w_0 < x < w_1$ and $t$ real. $j(s)$ is an "analytic characteristic function on V" when $j(s)$ is analytic on V and $j(iz)$ is positive definite in $z$ on the horizontal strip $-iV$, $j(s) = \int d(y)e^{sy}d(\mu)$ with $\mu$ a positive measure on the real line. Say $w_0 < x < w_1$. Then $j(x + it)/j(x)$ is the characteristic function relative to $t$ of the probability measure $\mu_x(S) := \int_S d(y)(e^{xy}/j(s))d(\mu)$. See E. Lukacs [25].

**Definition of a meromorphic characteristic function on C.**
Let us say that $h(z)$ is a "meromorphic characteristic function on C" when each of the following conditions hold.
(1) $h(z)$ is the reciprocal of an entire function, $j(z)$ say.

(2) All of the zeros of j(z) are on the union of the real axis and the imaginary axis.

(3) j(0) = 0. j(it) = 0 for at least one nonzero real t.

(4) The real zeros of j(z) are unbounded from above and also from below. Each real zero is simple.

(5) Let the successive real zeros be $w_k$ with k any integer, $w_0 = 0$ and $w_k < w_{k+1}$. Either h(z) or –h(z) is an analytic characteristic function on the open strip bounded by the vertical lines through $w_k$, $w_{k+1}$.

The Conjecture 5, C5, introduced below coupled with C^ serves to obtain that f(s) is an analytic characteristic function on $V_0$. The additional part that also assuming C5 brings is $g_0(y) > 0$ for y < 0.

The unconditional **Theorem 3.3** of **§3** implies the following. Given any r with $0 < r < P_0(\pi)$ there is a unique θ(r) with $0 < θ(r) < \pi$ and $P_0(θ(r)) = r$.

**C5 = Conjecture 5**.

(i) $0 < -A + c(0)$

(ii) $-A + c(0) < P_0(\pi)$.

**Definition of $v_0$.** Assume (i) and (ii) Set $v_0 := (1/\pi)θ(-A + c(0))$.

(iii) $P_0(v) > 0$, if $\pi < v \le \pi/v_0$.

The Conjecture 5 is corroborated by computer calculations done by the author using tables of certain of the $\gamma_k$ prepared by A. Odlyzko [33].

Motivation for C5 is presented in Part I, §7.

**§5 Context and statement of the main conditional theorem**

**Introduction**

The main conditional theorem assumes C^ and determines the representation $f(s) = \int_R d(y)e^{sy}g_0(y)$ valid on $V_0$, referred to in the introduction to §4. The Main conditional theorem implies RH and SZC.

A heuristic derivation of a formal expression for $g_0(y)$ is adumbrated next. In the formal partial fraction expansion f(s) = p(s) of (4.3), replace each term $1/(s – z)$ with its Laplace representation. Those representations are obtained from

$$1/s = \int_{y<0} (dy)e^{sy}, \text{ for Re(s) > 0, and } 1/s = \int_{y>0} (dy)e^{sy}(-1), \text{ for Re(s) < 0,}$$

(which are equivalent).

The heuristic result for y > 0 is $g_0(y) = P_0(\pi e^{-2y})$. When y < 0, it involves a formal series: $g_0(y) = (\sum_{k \geq 1} c(i\gamma_k)\cos(\gamma_k y)) + c(0) - P_0(\pi e^{2y})$.

**Definitions of λ(y), $g_0$(y).** Assume A is finite. $\sum_{k \geq 1} |c(i\gamma_k)\cos(\gamma_k y)|$ converges uniformly on the real line. Hence so does $\sum_{k \geq 1} c(i\gamma_k)\cos(\gamma_k y)$. (See §4, (4.2).)

$$\lambda(y) := 2\sum_{k \geq 1} c(i\gamma_k)\cos(\gamma_k y), \text{ for real y.}$$

$$g_0(y) := P_0(\pi e^{-2y}), \text{ for } y > 0. \quad g_0(y) := \lambda(y) + c(0) - P_0(\pi e^{2y}), \text{ for } y < 0.$$

λ(y) is continuous. $-A \leq \lambda(y) \leq A$. $\lambda(-y) = \lambda(y)$.

**Conditional claim 5.1** *Assume A is finite. Let Re(s) > 0.*

$$p_i(s) = \int_{y < 0} d(y)e^{sy}\lambda(y) = (-1)\int_{y > 0} d(y)\sinh(sy)\lambda(y).$$

**Proof** The method just described yields $p_i(s) = \sum_{k \geq 1} \int_{y < 0} d(y)t(k, s, y)$, with t(k, s, y) := $e^{sy}c(i\gamma_k)\cos(\gamma_k y)$. Now $|t(k, s, y)| \leq |c(i\gamma_k)| \cdot e^{xy}$, with x := Re(s). x > 0. A is finite. So $\sum_{k \geq 1} \int_{y < 0} d(y)|t(k, s, y)|$ is finite. Then $\sum_{k \geq 1} \int_{y < 0} t(k, s, y) = \int_{y < 0} \sum_{k \geq 1} t(k, s, y)$.

**Conditional corollary 5.1** *Assume $A < \infty$. On $V_0$: $p(s) = \int_R d(y)e^{sy}g_0(y)$.*

**Proof** $p(s) = p_i(s) + p_r(s)$. Together $0 < Re(s) < 4$ and the Corollary 3.2 of §3, give the representation of $p_r(s)$. That of $p_i(s)$ issues from $A < \infty$ and the previous Conditional claim 5.1.

In Part III, §2, Conditional theorem 2.2 and Part IV, §1, Conditional theorem 1.1 the respective culmination of the proof of each of the following theorems is achieved.

**Conditional theorem 5.1 Partial fraction representation of f(s).** *Assume $C^\wedge$. f(s) = p(s) on $C - Z°$.*

**Conditional theorem 5.2** *Assume $C^\wedge$. On $V_0$: $f(s) = \int_R d(y)e^{sy}g_0(y)$.*

**Definition of h(y).** Assume A is finite. Let y be real. Set h(y) := $\lambda(y) + c(0) - P_0(\pi e^{2y})$.

h(y) = $g_0(y)$, when y < 0. If $y \leq 0$, then $A + c(0) \geq h(y)$. Also $\liminf_{y \to -\infty} h(y) \geq -A + c(0)$.

If for any ε > 0, $|\zeta'(½ + it)| < e^{\varepsilon t}$ for all large positive t, then for nonreal z both (1)

and (2) below hold.
(1) Give any q with $0 < q < |Im(z)|$, one has $|c(i\gamma_k)\cos(\gamma_k z)| > \exp(\gamma_k q)$, for all large k.
(2) $\sum_{k \geq 1} c(i\gamma_k)\cos(\gamma_k z)$ diverges, since $\gamma_k \to \infty$.

The formal series $2\sum_{k \geq 1} c(i\gamma_k)(\cos(\gamma_k y))'$ for $\lambda'(y)$ is not absolutely convergent, if $|\zeta'(\frac{1}{2} + it)| = O(t^{1/4})$ for all large positive t.

Part (i) of C5 implies that A is finite.

The above analysis of the convergence and divergence properties of the series $\sum_{k \geq 1} c(i\gamma_k)\cos(\gamma_k z)$ also reveals similar behavior for $\sum_{k \geq 1} c(i\gamma_k)\sin(\gamma_k z)$. The same holds for $\sum_{k \geq 1} c(i\gamma_k)\exp(i\gamma_k z)$, at z with $Im(z) < 0$.

**Definition of e(z).** If $A < \infty$, then $\sum_{k \geq 1} c(i\gamma_k)\exp(\gamma_k z)$, for z with $Re(z) \leq 0$, converges uniformly to a continuous function e(z).

Say $\delta \geq 0$. Let $H_\delta$ be the half-plane of z with $Re(z) \leq -\delta$. At each z in $H_\delta$:

$$|e(z)| \leq \sum_{k \geq 1} |c(i\gamma_k)\exp(\gamma_k z)| \leq \sum_{k \geq 1} |c(i\gamma_k)|\exp(-\gamma_k\delta) \leq (A.2)\exp(-\gamma_1\delta).$$

C2 implies the C2′ stated next.

C2′ : $\sum_{k \geq 1} |c(i\gamma_k)|\exp(-\gamma_k\varepsilon)$ is finite, for arbitrarily small positive $\varepsilon$.

C2′ implies that $\sum_{k \geq 1} |c(i\gamma_k)\exp(\gamma_k z)|$ converges for $Re(z) < 0$ and that the convergence is uniform on $H_\delta$ when $\delta > 0$. Then e(z) is analytic for the z with $Re(z) < 0$. Also e(z) vanishes rapidly uniformly, as $Re(z) \to -\infty$.

Assume C2. Say y is real. Set $\lambda_S(y) := \sum_{k \geq 1} c(i\gamma_k)\sin(\gamma_k y)$ and $\lambda_C(y) := \lambda(y)/2$. One has $e(iy) = \lambda_C(y) + i\lambda_S(y)$ and $\lambda_S(-y) = -\lambda_S(y)$.

Constituent themes of the main conditional theorem are treated next. Each theme is initiated by restating relevant unconditional results. Next complimentary conditional results are given. Then the issue of the synthesis of those unconditional and conditional results is presented.

**(5.1) Analyticity of f(s)**

(i) **Unconditional result.** f(s) is analytic on an open set containing all of the s for which each of $|Re(s)| \geq \frac{1}{2}$ and $s \neq 4w$ for $w \geq 1$ holds.

(ii) **Conditional result.** Assume RH. f(s) is analytic for s with |Re(s)| < ½ and Re(s) ≠ 0.

(iii) **Synthesis.** RH implies f(s) is analytic at s if Re(s) ≠ 0 and s is not a multiple of four.

## (5.2) Integrability of |f(s)| on vertical lines

(i) **Unconditional result.** |f(s)| is integrable and square integrable on any vertical line x + i**R**, with x real, |x| ≥ ½ and x ≠ ±4w for w ≥ 1.

(ii) **Conditional result.** Assume RH and C1. Then |f(s)| is integrable and square integrable on any vertical line x + i**R**, with -½ < x < ½, other than the imaginary axis.

(iii) **Synthesis.** Assume RH and C1. Then |f(s)| is in $L_1 \cap L_2$ on x + i**R**, if x is real, but not a multiple of four.

## (5.3) The Laplace densities g, $g_0$.

We seek to determine the densities g, $g_0$ for which the following Laplace representations hold. On $V_0'$: $f(s) = \int_R d(y) e^{sy} g(y)$. On $V_0$: $f(s) = \int_R d(y) e^{sy} g_0(y)$.

(i) **Unconditional results.** g(y) is the restriction to the real line of the entire function $P_0(\pi e^{-2z})$. $|P_0(\pi e^{-2z})| \sim O(\exp(-4z_1))$, with $z_1 := Re(z)$ and $z_1 \to \infty$. Therefore $\int_{y>0} d(y) |e^{sy} P_0(\pi e^{-2y})|$ is finite, if Re(s) < 4. When x > ½, $e^{xy}|P_0(\pi e^{-2y})|$ vanishes rapidly on the real line as y → -∞. (See §3, Theorem 3.2 (ii).) Therefore $\int_{y<0} d(y) e^{xy} |P_0(\pi e^{-2y})|$ is finite.

**Review** Introduction of §5.

(ii) **Conditional results for $g_0(y)$.** Assume C2. $g_0(y) = P_0(\pi e^{-2y})$, on the ray y > 0. As in (i) above: $\int_{y>0} d(y) |e^{sy} P_0(\pi e^{-2y})|$ is finite, if Re(s) < 4. λ(y) is defined and continuous for real y. Assume y < 0. Then $g_0(y) := \lambda(y) + c(0) - P_0(\pi e^{-2y})$. Thus $g_0(y)$ is continuous and bounded for y < 0. Hence for x > 0, $e^{xy} g_0(y)$ decays rapidly as y → -∞. Thus $\int_{y<0} d(y) |e^{sy} g_0(y)|$ is finite, if Re(s) > 0.

(iii) **Synthesis.** Assume C2. Then $\int_R d(y) |e^{sy} g_0(y)|$ is finite for s on $V_0$, since the integrand vanishes rapidly as y → ±∞.

## (5.4) The representation of f(s) as a Laplace transform

(i) **The main unconditional theorem.** On the strip $V_0'$ of s with ½ < Re(s) < 4: $f(s) = \int_{\mathbf{R}} d(y) e^{sy} P_0(\pi e^{-2y})$. $P_0(v)$ is bounded, for $0 \le v \le \pi$. See the Main unconditional theorem (1) of §3 and its proof in Results when β = ¼, §6, Part II.

(ii) **The conditional representation.** Assume C^. On the strip $V_0$ of s with 0 < Re(s) < 4: $f(s) = \int_{\mathbf{R}} d(y) e^{sy} g_0(y)$. Here $g_0(y) = P_0(\pi e^{-2y})$, if y > 0. If y < 0, then $g_0(y) := \lambda(y) + c(0) - P_0(\pi e^{2y})$. $g_0(y)$ is bounded on the nonzero reals. See the Conditional theorem 5.2 of §5 and its proof as Conditional theorem 1.1 in Part IV, §1.

(iii) **Synthesis.** Combining (i) and (ii) for s in $V_0'$ gives the following new results.

Assume C^. On $V_0'$, $f(s) - F(s) = \int_{y<0} d(y) e^{sy} g_0(y) = \int_{y<0} d(y) e^{sy} P_0(\pi e^{-2y})$. Thus $g_0(y) = P_0(\pi e^{-2y})$ also for y < 0. Next apply the conditional representation on all of $V_0$ to obtain the following.

**Main conditional theorem (1)** *Assume C^.*
**(i)  The equality of the conditional and unconditional Laplace densities.**
*If y is real, then:*

<div align="right">**Eq. (\*)**</div>

$$\lambda(y) + c(0) - P_0(\pi e^{2y}) = P_0(\pi e^{-2y}).$$

**(i′) The boundedness of the density.**
*$P_0(v)$ is bounded on the real axis.*
**(ii) The conditional extension of the unconditional Laplace representation of f(s) on $V_0'$ to $V_0$.**
*On the strip $V_0$: $f(s) = \int_{\mathbf{R}} d(y) e^{sy} P_0(\pi e^{-2y})$.*

See Part IV, §2 Proof of the Main conditional theorem (1).

**Conditional criterion for continuity of $g_0$ at the origin.** Assume C^. It follows that $\lim_{y<0, y \to 0} g_0(y) = \lim_{y>0, y \to 0} g_0(y)$ holds when:

$$\sum_{k \ge 1} c(i\gamma_k) = -(c(0)/2 + \sum_{k \ge 1} c(4k)).$$

The previous criterion and the Conditional theorem 5.3 (1) stated next follow from the Main conditional theorem (1).

**Conditional theorem 5.3**
(1) *Assume C^. $g_0(z) = P_0(\pi e^{-2z})$, respectively*

$$\lambda(z) = -c(0) + P_0(\pi e^{2z}) + P_0(\pi e^{-2z}) = -(c(0) + 2\sum_{k \geq 1} c(4k)cosh(4ky)),$$

*holds on the real line and so extends $g_0$, respectively $\lambda$, to an entire function on C of period $i\pi/2$.*
(2) *Assume RH and A is finite. If $|Im(s)| < \gamma_1$, then $p_i(s) = \int_{y>0} d(y)sin(sy)2e(-y)$.*

**Proof of (2)** Say s = iz. $p_i(iz) = \sum_{k \geq 1} \int_{y>0} d(y)\theta(k, s, y)$, with $\theta(k, s, y) :=$ $c(i\gamma_k)2sin(sy)exp(-\gamma_k y)$. Set t := Im(s). Say $0 \leq \delta < \gamma_1$ and $|t| \leq \delta$. Now $|\theta(k, s, y)|$ $\leq 2|c(i\gamma_k)|\cdot exp(-(\gamma_1 - |t|)y)$. Also $\gamma_1 - |t| \geq u$, with $u = \gamma_1 - \delta$. A is finite. So $\sum_{k \geq 1} \int_{y>0} d(y)|\theta(k, s, y)| \leq A\cdot\int_{y>0} d(y)e^{-uy} < \infty$. So $\sum_{k \geq 1} \int_{y>0} \theta(k, s, y) = \int_{y>0} \sum_{k \geq 1} \theta(k, s, y)$.

**§6 Relations of $\gamma_n$, $\zeta'(\frac{1}{2} + i\gamma_n)$ to their predecessors and the $\zeta(\frac{1}{2} + 4k)$.**

The Conditional lemma 4.1 of Part I, §4, (4.2), established that C2 implies A < ∞.

**Review** §5, Introduction, Definition of e(z).

**Conditional lemma 6.1** *Assume $A := 2\sum_{k \geq 1} |c(i\gamma_k)|$ is finite. Say Re(z) > 0.*

**Eq. (°)**
$$e(-z) = (z/\pi)\int_{y>0} d(y)(1/(z^2 + y^2))\lambda(y), \text{ with } \lambda(y) := 2\sum_{k \geq 1} c(i\gamma_k)cos(\gamma_k y).$$

**Proof** Say Re(z) > 0 and $j \geq 0$. $e^{-zj} = (2z/\pi)\int_{y>0} d(y)(1/(z^2 + y^2))cos(jy)$. Assume C2. Fix z.
$e(-z) = \sum_{k \geq 1} \int_{y>0} d(y)t(k, y)$, with $t(k, y) := (2/\pi)(z/(z^2 + y^2)c(i\gamma_k)cos(\gamma_k y)$.

Claim $\sum_{k \geq 1} \int_{y>0} d(y)|t(k, y)| < \infty$. $\sum_{k \geq 1} \int_{y>0} d(y)\theta(k, y) = \int_{y>0} d(y)\sum_{k \geq 1} \theta(k, y)$, *for $\theta = |t|$ and thus for $\theta = t$.*

Proof of Claim $\sum_{k \geq 1} \int_{y>0} d(y)|t(k, y)| \leq (A/\pi)|z|\int_{y>0} d(y)/|z^2 + y^2|$. Now $z = |z|\cdot e^{i\varphi}$, with $|z| \neq 0$ and $|\varphi| < \pi/2$. Set $u = e^{i\cdot|\varphi|}$.

Subclaim *$(1/\pi)|z| \int_{y>0} d(y)/|z^2 + y^2| \leq (1/cos(|\varphi|))(\frac{1}{2} + |\varphi|/\pi).$*

Proof of Subclaim $|z|\int_{y>0} d(y)/|z^2 + y^2| = \int_{y>0} d(y)/(|u - iy|\cdot|u - (-iy)|) \leq \int_{y>0} d(y)|u - iy|^{-2} = (1/cos(|\varphi|))\int_{y>-tan(|\varphi|)} d(y)/(1 + y^2) = (1/cos(|\varphi|))(\pi/2 + |\varphi|)$.

**Definition of j(u).** Set $j(u) := -c(0) + P_0(\pi e^{2u}) + P_0(\pi e^{-2u})$, for complex u.

j(u) is an entire function.

**Review** §3, Theorem 3.2 (ii). §4, Conditional corollary 4.1. The definition of C^ in §4, (4.4). The Main conditional theorem (1) (i′) of §5, (5.4).

The following Conditional corollary 6.1 is obtained in Part V, §1, as Conditional corollary 1.

**Conditional corollary 6.1** *Assume*
*(′): $\theta < 1$ and for $v \geq 0$, $|P_0(v)| \sim O((log(v))^\theta)$, as $v \to \infty$.*
*Then:*
*(1) $(z/\pi)\int_{y \geq 0} d(y)(1/(y^2 + z^2))j(y)$ converges absolutely to an analytic function on the half-plane $Re(z) > 0$.*

**Definition of υ(z).** Assume the previous corollary, (′) thereof and $Re(z) > 0$. Set

$$\upsilon(z) := (z/\pi)\int_{y > 0} d(y)(1/(z^2 + y^2))(-c(0) + P_0(\pi e^{2y}) + P_0(\pi e^{-2y})).$$

**A key observation** In the previous integrand the only s for which $\zeta(s)$ enters therein are ½, from the critical strip, and the ½ + 4k with k a positive integer. See Part I, §3, for the comment (2) on the Taylor series for $P_0(z)$.

Attempting to simplify the formula used to define υ(z) by a term by term integration using the Taylor series for $P_0(v)$ encounters the obstacle, posed by $P_0(\pi e^{2y})$, of the exponential rates of divergence to infinity as $y \to \infty$ of the $(e^{2y})^{2k}$ arising from the positive integers k.

**Conditional corollary 6.1** *Assume (′) above.*
*(2) υ(z) extends to an entire function on C.*

Review Eq. (*) of §5, (5.4), Main conditional theorem (1) (i).

**Conditional claim 6.1** *Assume the validity of Conditional corollary 6.1 and of*
*(*): A is finite and for $y > 0$, Eq. (*) holds, $\lambda(y) = j(y)$.*
*(1) $e(-z) = \upsilon(z)$ on the half-plane $Re(z) > 0$:*

$$\sum_{k \geq 1} c(i\gamma_k)exp(-\gamma_k z) = (z/\pi)\int_{y > 0} d(y)(1/(z^2 + y^2))(-c(0) + P_0(\pi e^{2y}) + P_0(\pi e^{-2y}),$$

*with the sum and the integral absolutely convergent.*

(2) *e(s) has an analytic extension to the entire complex plane.*

**Proof** Assume A is finite. Conditional lemma 6.1 yields Eq. (°). Also assume Eq. (*) holds. In the integrand of Eq. (°) replace $\lambda(y)$ with $j(y)$. $\lambda(y)$ is bounded for $y > 0$. The Conditional claim 6.1 then follows from the Conditional corollary 6.1.

(1), (2) are derived from (*) in Conditional corollary 2 of Part V, §1.

**Conditional claim 6.2** *Assume the validity of C^, of the Main conditional theorem (1) (i) of §5, (5.4), and of Conditional corollary 6.1. Then (1), (2) of Conditional claim 6.1 hold.*

**Proof** Assume C^. Then $A < \infty$. (See (4.4).) The Main conditional theorem (1) (i) then gives $\lambda(y) = j(y)$, for real y. Assume Conditional corollary 6.1. (1), (2) follow.

The Main conditional theorem (1) (i) is established in Part IV, §2. (1), (2) of Conditional claim 6.1 are obtained from C^ in Conditional corollary 3 (1) of Part V, §1.

**Definitions of e(z, n), ê(z, n).** Assume A is finite. Set

$$e(z, n) := \sum_{1 \leq k \leq n} c(i\gamma_k)\exp(\gamma_k z).$$

Say $\mathrm{Re}(z) \leq 0$. Let

$$\hat{e}(z, n) := \sum_{k \geq n} c(i\gamma_k)\exp(\gamma_k z) = e(z) - e(z, n - 1).$$

If RH and SZC hold, then $c(i\gamma_k) = 1/(b(i\gamma_k)\zeta'(\tfrac{1}{2} + i\gamma_k))$. Assume RH and that A is finite. The $\gamma_n$, $\zeta'(\tfrac{1}{2} + i\gamma_n)$ of the sequence $\gamma_1$, $\zeta'(\tfrac{1}{2} + i\gamma_1)$; ...$\gamma_n$, $\zeta'(\tfrac{1}{2} + i\gamma_n)$;... can successively be expressed in terms of their predecessors $\gamma_k$, $\zeta'(\tfrac{1}{2} + i\gamma_k)$, with $1 \leq k \leq n - 1$, and, in the case of $\zeta'(\tfrac{1}{2} + i\gamma_n)$, also in terms of $\gamma_n$. That involves e(z) and limit processes. Say $x > 0$. $\gamma_1 = -\lim_{x \to \infty} (1/x)\log(-e(-x))$. $\zeta'(\tfrac{1}{2} + i\gamma_1) = 1/(b(i\gamma_1)\lim_{\mathrm{Re}(z) \to \infty} \exp(\gamma_1 z) \cdot e(-z))$. $\gamma_n = -\lim_{x \to \infty} (1/x)\log((-1)^n \cdot \hat{e}(-x, n))$. Also $\zeta'(\tfrac{1}{2} + i\gamma_n) = 1/(b(i\gamma_n)\lim_{\mathrm{Re}(z) \to \infty} \exp(\gamma_n z)\hat{e}(-z, n))$. In those expressions e(-z) can be replaced with $\upsilon(z)$, when $\mathrm{Re}(z) > 0$, provided the previous Conditional claim 6.1 (1) and Eq. (*) thereof hold. Note that the $\zeta(\tfrac{1}{2} + 4k)$, with k a nonnegative integer, enter into $\upsilon(z)$ via the coefficients in the Taylor series for $P_0(v)$. The results thereby obtained are stated next.

**Conditional corollary 6.2** *Assume the validity of RH, Conditional corollary 6.1 (1) and (\*) of Conditional claim 6.1 .*

**Relations of $\gamma_n$, $\zeta'(\frac{1}{2} + i\gamma_n)$ to their predecessors and the $\zeta(\frac{1}{2} + 4k)$**

$$\gamma_1 = -lim_{x>0,\, x\to\infty} (1/x)log((-1)^n \cdot v(x)).$$

$$\zeta'(\frac{1}{2} + i\gamma_1) = 1/(b(i\gamma_1)lim_{Re(z)\to\infty} exp(\gamma_1 z)\cdot v(z)).$$

$$\gamma_n = -lim_{x>0,\, x\to\infty} (1/x)log((-1)^n (v(x) - e(-x, n-1))).$$

$$\zeta'(\frac{1}{2} + i\gamma_n) = 1/(b(i\gamma_n)lim_{Re(z)\to\infty} exp(\gamma_n z)(v(z) - e(-z, n-1))).$$

The previous relations are established using C^ in Conditional corollary 3 (2) of Part V, §1. That corollary is stated next.

**Conditional corollary 6.3** *Assume C^.*
**(1)** *$e(-z) = v(z)$, provided $Re(z) > 0$. $v(z)$ extends $e(-z)$ to an entire function on C.*
**(2)** *The above relations hold for $\gamma_n$ and $\zeta'(\frac{1}{2} + i\gamma_n)$.*

**Review** §3 Lemma 3.1. See Conditional claim 4.1 of §4, (4.3).

Assume A is finite. Say s is distinct from the $i\gamma_k$, with $k \geq 1$. $\sum_{k\geq1} |c(i\gamma_k)/(s - i\gamma_k)|$ converges.
**Definition of $p_{i,+}(s)$.**
$$p_{i,+}(s) := \sum_{k\geq1} c(i\gamma_k)/(s - i\gamma_k).$$

$p_{i,+}(s)$ is analytic except for simple poles at $i\gamma_k$ for $k \geq 1$. Assume RH and that A is finite. Then $p_i(s) = p_{i,+}(s) - p_{i,+}(-s)$.

$\sum_{k\geq1} \gamma_k^{-p}$ is finite, for $p > 1$. $\xi(\frac{1}{2}) > 0$.

**Definition of $\Xi(s)$.** Set $\Xi(s) := (\xi(\frac{1}{2}))^{\frac{1}{2}} \cdot \Pi_{k\geq1} (1 - (s/(i\gamma_k)))exp(s/(i\gamma_k))$.

$p_{i,+}(s) = \Lambda(s)/\Xi(s)$, with $\Lambda(s)$ entire.

Assume RH and SZC for the following factorization and fission by symmetry.

**The Hadamard factorization of $\xi(\frac{1}{2} + s)$.**

$$\xi(\frac{1}{2} + s) = \xi(\frac{1}{2}) \cdot \Pi_{k\geq1} (1 + (s/\gamma_k)^2).$$

The previous factorization follows from Theorem 9.27 on p 205 of G. Everest, T. Ward [40].

**ξ-fission**

$$\xi(½ + s) = \mathbf{\Xi}(s)\mathbf{\Xi}(-s).$$

**Conditional corollary 6.4  Laplace transform representation of $ip_{i,+}$(is).**
*If A is finite and Re(s) < $\gamma_l$, then $ip_{i,+}$(is) = $\int_{y > 0} d(y)e^{sy} e(-y)$.*

**Proof** Say x := Re(s) < $\gamma_1$. Then $\sum_{k \geq 1} \int_{y > 0} d(y)|c(i\gamma_k)|\exp((x - \gamma_k)y) \leq A/(\gamma_1\text{-x})$. Also assume A is finite. Then in $\sum_{k \geq 1} \int_{y > 0} d(y)c(i\gamma_k)\exp((s- \gamma_k)y)$ apply the interchange $\sum_{k \geq 1} \int_{y > 0} d(y) = \int_{y > 0} d(y)\sum_{k \geq 1}$.

Consider the principal branch of z with z nonzero and |arg(z)| < π. Let q(t) be a path which begins at z, never crosses the negative real axis, does not pass through the origin and has Re(q( y)) → ∞ as y → ∞. $E_1(z) := \int_q d(y)e^{-y}$ /y defines the exponential integral. If Re(z) > 0, then $E_1(z) := \int_{y \geq 1} d(y)e^{-zy}$ /y.

**Definition of $\Theta(\theta, z)$.** Assume Im(z) < 0 and θ > 0.

$$\Theta(\theta, z) := (1/(2i))\cdot\sum_{\sigma = \pm 1} e^{-\sigma\theta z}\cdot E_1(-\sigma\theta z)$$

Review Conditional claim 6.1.

**Conditional corollary 6.5 Representation of $p_{i,+}$(z) via j(y).**
(1) *Assume A is finite and for y > 0, λ(y) = j(y). Say Im(z) < 0. Then*

$$p_{i,+} (z) = \int_{\theta > 0} d(\theta)j(\theta)\Theta(\theta, z).$$

(2) *Assume $\hat{C}$. Then the previous representation of $p_{i,+}$(z) holds on the lower half-plane of z with Im(z) < 0.*

Conditional corollary 6.5 is proven as Conditional corollary 4 in §2 of Part V. The proof  of (2) relies on Main conditional theorem (1) (i), proven in Part IV, §2.

**§7 f(s) as a meromorphic characteristic function on C**

(i) **Unconditional result.** $P_0$(v) > 0, when 0 < v ≤ π. by Theorem 3.3.

We provide motivation for C5 by next establishing the Conditional lemma 7.2 which leads from Main conditional theorem (1) (i) of §5, (5.4), to Main conditional theorem (2) below.

**Review** Definition of h(y) of §5, Introduction. Eq. (*) of §5, (5.4), Main conditional theorem (1) (i). (See Eq. (*) in §6, Conditional corollary 6.2.)

Assume A is finite. $\liminf_{y \to -\infty} h(y) \geq$ -A + c(0).

Assume C^. Apply Main conditional theorem (1) (i). One obtains Eq. (*): $h(y) = P_0(\pi e^{-2y})$, for real y.

It results from Eq. (*) that $-A \leq -c(0) + P_0(\pi/v) + P_0(\pi v) \leq A$, when v > 0. In particular

**Eq. (')**

$$P_0(\pi v) \geq c(0) - A - P_0(\pi/v).$$

Apply Theorem 3.3 to obtain that for $v \geq 1$, $c(0) - A - P_0(\pi/v)$ is strictly increasing from $c(0) - A - P_0(\pi)$ towards $c(0) - A$. If $c(0) - A - P_0(\pi) > 0$, then $P_0(v) > 0$, when $v \geq \pi$. Calculations support C5 of §4, (4.5).

**Claim 7.1** *Suppose that the $\gamma_k$ with $k \geq 1$ are independent over the rationals (integers). Then for any real T the values of $\lambda(y)$ arising from the $y < -T$ are dense in [-A, A]. Hence $\liminf_{y \to -\infty} h(y) = -A + c(0)$.*

The claim follows from an approximation theorem of Kronecker. That theorem is presented in Apostle [4].

Assume (i) and (ii) of C5. The $v_0$ of C5 satisfies $0 < v_0 < 1$ and $c(0) - A - P_0(\pi v_0) = 0$. Thus $c(0) - A - P_0(\pi/v) > 0$, for $v > 1/v_0$. Also assume that $h(y) = P_0(\pi e^{-2y})$, for real y. Apply Eq. (') to obtain that $P_0(v) > 0$, for $v > \pi/v_0$. Also assume (iii) of C5. Then $P_0(v) > 0$, if $\pi < v \leq \pi/v_0$.

(ii) **Conditional result on positivity.**

**Conditional lemma 7.1** *Assume C5 and Eq. (*) hold. Then $\inf_{v > \pi} P_0(v) > 0$ and $\liminf_{v \to \infty} P_0(v) \geq -A + c(0) > 0$.*

(iii) **Synthesis. Positivity of $P_0(v)$.**

**Conditional lemma 7.2** *Assume C5 and Eq. (\*) hold. Then $P_0(v) > 0$, for all $v > 0$. Say $\varepsilon > 0$. Then $\inf_{v > \varepsilon} P_0(v) > 0$.*

**The Main conditional theorems (2)-(3).**

Assume the validity the Main conditional theorem (1) (i). The following theorem results via Conditional lemma 7.2.

**Main conditional theorem**
**(2) *Assume C^ and C5. f(s) is an analytic characteristic function on $V_0$:***

$$f(s) := 1/n(s) = \int_R d(y)e^{sy}P_0(\pi e^{-2y}),$$

**with $P_0(v)$ positive for $v > 0$. Also $\inf_{v > \varepsilon} P_0(v) > 0$, for any $\varepsilon > 0$.**

(i) **Unconditional results.** Say $w \geq 1$ and $s$ is in $V_{4w}$. Then the analytic function $(-1)^w f(s)$ has a Laplace transform representation with positive density $P_{4w}(\pi e^{-2y})$. Thus $(-1)^w f(s)$ is an analytic characteristic function on $V_{4w}$.

Applying the relation $f(-s) = -f(s)$ one obtains the counterparts of the above results for the negative half-plane of $s$ with $\text{Re}(s) < 0$. When $w \leq -1$, set $P_{4w}(\pi v) := P_{-4(w+1)}(\pi/v)$, for $v > 0$.

The conjunction of the previous Main conditional theorem (2) and the Main unconditional theorem (4) of Part I, §3, engenders the following.

**Main conditional theorem**
**(3) *Assume C^ and C5. f(s) is a meromorphic characteristic function on the complex plane: When w is an integer and s is in $V_{4w}$,***

$$(-1)^w f(s) = \int_R d(y)e^{sy}P_{4w}(\pi e^{-2y}),$$

**with $P_{4w}(z)$ entire in $z$ and $P_{4w}(v)$ positive for $v > 0$.**

See Part IV, §3 Proofs of the Main conditional theorems (2)-(3).

**Definitions of ridge, groove functions.** Say $q$ maps the subset $V$ of the complex plane into $C$ extended by a point at infinity. $Q$ is a "ridge function" on $V$ (with ridge contained in the real axis) shall mean that (i), (ii) as follows hold.
(i) If $x$, $t$ are nonzero reals and $x$, $x + it$ are in $V$, then $|q(x + it)| < |q(x)|$.

(ii) If t is a nonzero real and each of 0, i·t is in V, then ⎮q(i·t)⎮ ≤ ⎮q(0)⎮. (See E. Lukacs [25], Ch. 7, p. 195). Reversing the prior inequalities defines when q has the "groove property" on V.

The "extended sense" of the ridge property allows ⎮q(x + it)⎮ ≤ ⎮q(x)⎮. That for the groove property allows ⎮q(x + it)⎮ ≥ ⎮q(x)⎮.

Each meromorphic characteristic function on the complex plane is also a ridge function on C.

The entire function $n(s)$ has a simple zero at 0 and does not vanish at the x with $0 < x < 4$.

**The Main conditional theorem (2) implies RH and SZC.**

The Main conditional theorem (2) entails that $f(s)$ is a ridge function on $V_0$. Thus $n(s)$ is a groove function on $V_0$. Therefore $n(s) \neq 0$ on $V_0$. Hence RH holds. SZC also holds, as is explained next.

**Claim 7.2** *Let V be an open subset of C. Say j(z) maps V into C, $z_0$ is a zero of j on V and $x_0 := Re(z_0)$ is also in V. Assume j(z) is an extended groove function on the open subset V′ of z in V with Re(z) > $x_0$.*
*(i) If j(z) is continuous (from the right) at $x_0$, $z_0$ then $j(x_0) = 0$.*
*(ii) If j(z) is analytic on neighborhoods in V of $x_0$, $z_0$ and $j′(x_0)$ is nonzero, then $j′(z_0) \neq 0$.*

**Proof of Claim 7.2 (ii)**. $j(x_0) = 0$. Also, for small h > 0, the approximations $(1/h)j(\theta + h) = j′(\theta) + hO(h)$ hold for $\theta = x_0$, $z_0$. The extended groove property gives ⎮j($z_0$ + h)⎮ ≥ ⎮j($x_0$ + h)⎮. Thus ⎮j′($z_0$)⎮ ≥ ⎮j′($x_0$)⎮. Then $j′(x_0) \neq 0$ yields $j′(z_0) \neq 0$.

**The role of b(s) in the groove property of n(s).** Say t is real. If x > 1 and t is nonzero then $\zeta(x) > ⎮\zeta(x + it)⎮$. Computer calculations indicate that there exist x with ½ < x < 1 and nonzero real $t_0$ for which:

$$\zeta(x + it_0) = \zeta(x) \text{ and } (d/dt)(⎮\zeta(x + it)⎮^2) \text{ is nonzero at } t_0.$$

Thus $\zeta(s)$ is neither a ridge nor a groove function on the strip V(½, 1). This underscores the significance of the factor b(s) in obtaining (from C^ and C5) the groove property of $n(s) := b(s)\zeta(½ + s)$ on $V_0$.

**Polynomial counterexamples for the groove property**.

Next we specify polynomials $n_1(s)$ which share properties in common with $n(s)$, but do not have the corresponding groove property.

Let $\alpha(s) := \Pi_{1 \leq w \leq M}(s - j_w)$, $\beta(s) := \Pi_{1 \leq k \leq N}(s - v_k)$, with the zeros positive and increasing with the index w, k respectively. Set $p(s) := s\alpha(s)\alpha(-s)$, $q(s) := \beta(is)\beta(-is)$, $n_1(s) := p(s)q(s)$ and $f_1(s) := 1/n_1(s)$.

The roots of $n_1(s)$ are on the union of the real axis and the imaginary axis. Each root is simple. $n_1(0) = 0$. q is even in s. Also $p, n_1, f_1$ are odd. $n_1{}'$ is even. When $\theta = q, p, n_1, f_1$ or $n_1{}'$, one has $\theta(s^*) = (\theta(s))^*$. Also $n_1{}'$ is real-valued on the aforementioned axes.

$n_1(s) := p(s)q(s)$ may not have the groove property, even though $p(s)$ does. The factor q of $n_1$ can cause problems. For all small $x > 0$, $|n_1(x + iv)| < |n_1(x)|$ at each root iv of q for which $|n_1{}'(iv)| / n_1{}'(0) < 1$, since as x vanishes $n_1(x + iv)/n_1(x)$ approaches $n_1{}'(iv)/n_1{}'(0)$. A simple instance has $n_1(s) := s\omega(s)\omega(-s)$, with $\omega(s) := (s - iv_1)(s - iv_2)$. Then $|n_1{}'(iv)| / n_1{}'(0) = 2(1 - (v_1/v_2)^2)$. Thus $n_1$ is not a groove function, if $(v_2 - v_1)/v_1 < 2^{1/2} - 1$.

**Metric norms and analytic characteristic functions.**

The property of being an analytic characteristic function finds expression in a geometric context. This involves the internal metric norms obtained by the author in Part VI (see A. Csizmazia [12]), as a generalization of the Fourier metric norms of J. von Neumann and I. J. Schoenberg [28]. See I. J. Schoenberg [37]. An example engendered by $f(s) := 1/n(s)$ is delineated next.

**Review** §3, A geometric consequence of the Main unconditional theorem (4), Definitions of $m_x(t)$, $d_x(t_1, t_2)$.

The next assertion is proven in Part VI. That proof uses Conditional corollary 4.1 of Part IV, §4.

*Assume C^ and C5. Fix x with 0 < x < 4.*
*$m_x(t)$ is a (finite-valued) metric norm in t on the real line.*
*$d_x(t_1, t_2)$ is a translation invariant metric in $t_1$, $t_2$ on the real line.*

Refer to Metric norms and analytic characteristic functions of Part II, §6 and Part IV, §4.

## §8 Meromorphic characteristic functions arising from L(s, χ) or r(s)

We have described a context relative to f(s) := 1/(b(s)ζ(½ + s))  being a meromorphic characteristic function on C. An analogous context is outlined next for the recondite case of certain functions built from Dirichlet L-functions, L(s, χ), instead of the zeta-function. The context deriving from the Ramanujan tau Dirichlet function, r(s), is then considered.

The possibility of success in the case of L(s, χ) is obscured by the following.

**Lemma 8.1** *Assume that: either (1) or (2) below holds.*
*(1) χ(n) is nonreal, for some n.*
*(2) There exist primes p, p' such that χ(p) = 1 and χ(p') = -1.*

*Then for any x > 1 and T > 0, there exist real numbers (integers) t, t' > T, such that:*

$$|L(x + it, \chi)| < |L(x, \chi)| < |L(x + it', \chi)|.$$

Lemma 8.1 is proven as Corollary 3 in the Appendix herein. That relies on the Euler factorization, when Re(s) > 1,

$$L(s, \chi) = \Pi_{p:\ prime > 0.}\ (1 - \chi(p) \cdot p^{-s})^{-1}$$

(see T. M. Apostol [3]) and a Kronecker approximation theorem (see Appendix and T. M. Apostol [4]). There is computational evidence for certain χ that x, $t_0$ exist with ½ < x < 1 and $t_0$ a nonzero real, for which L(x + $it_0$, χ) = L(x, χ) and (d/dt)(|L(x + it, χ)|$^2$) is nonzero at $t_0$. Thus L(s, χ)) is neither a groove nor a ridge function on the vertical strip of s with ½ < Re(s) <  1.

The construction of a meromorphic characteristic function on C from L(s, χ) has a special feature when it is permitted that some of the values of χ be nonreal. Then the fundamental building block is L(s, χ)·L(s, χ*). (Note that L(s, χ*) = L(s*, χ)*.)

### (8.1) Dirichlet L-functions, L(s, χ)

Assume χ is a non-principal character mod k, with k > 1. Set

$$q(s, \chi) := (k/\pi)^{s/2}\Gamma((s + \theta)/2)L(s, \chi).$$

Here θ = 0 , if χ(-1) = 1, whereas θ = 1 , if χ(-1) = -1.

The functional equation is

$$q(1 - s, \chi) = \omega(\chi)q(s, \chi^*).$$

Here $\omega(\chi) := jG(1, \chi)(k^{-1/2})$. This j is 1, if $\chi(-1) = 1$, whereas j is –i, if $\chi(-1) = -1$. Moreover $G(1, \chi)$ is the Gaussian sum $\sum_{1 \leq r \leq k} \chi(r)\exp(2\pi r/k)$. Note that $|\omega(\chi)| = 1$ and $\omega(\chi^*) = (\omega(\chi))^*$.

**Dirichlet L-function conjecture, DLFC.** Each nonreal zero of the L-function $L(s, \chi)$ has real part ½.

**Simple zeros conjecture, SZCD.** The nonreal zeros of $L(s, \chi)$ are all simple.

Set $f(s, -1) := \sin(\pi s/4)$ and $f(s, 1) := s\cdot\cos(\pi s/4)$. Define

$$v(s, \chi) := f(s, \chi(-1))q(\frac{1}{2} + s, \chi).$$

Then $v(s, \chi)$ is an entire function of s. Also

$$v(-s, \chi) = -\omega(\chi)v(s, \chi^*).$$

Assume $\chi$ is real-valued. Set $n(s, \chi) = v(s, \chi)$. The next conjecture restricted to $0 < x < \frac{1}{2}$ implies DLFC and SZCD hold for $\chi$.

**Conjecture D1.** If $\chi$ is real-valued, then $n(s, \chi)$ is a groove function in s on the complex plane.

Now allow $\chi$ to be complex-valued. Define

$$g(s, \chi) := v(s, \chi)v(s, \chi^*) \text{ and } h(s, \chi) := s\cdot\sin(\pi s/2)q(\frac{1}{2} + s, \chi)q(\frac{1}{2} + s, \chi^*).$$

Then g, h are even functions of s. Fix $n = g$ or $n = h$. The next conjecture restricted to $0 < x < \frac{1}{2}$ implies DLFC and SZCD hold for $\chi$.

**Conjecture D2.** $n(s, \chi)$ is a groove function in s on the complex plane.

Assume that $\chi(n)$ is not real for some n. Set

$$j(s, \chi) := f(2s, \chi(-1))q(\frac{1}{2} + s, \chi)q(\frac{1}{2} + s, \chi^*).$$

The possibility that $j(s, \chi)$ is a groove function in s is sensitive to any exceptional nearness of pairs z, z′ with either both of z, z′ or both of z*, z′ zeros of $L(\frac{1}{2} + s, \chi)$.

Computer calculations support conjectures D1 and D2. Those conjectures are corollaries of the one which follows.

**Main conjecture for L(s, χ).** Let s = u, u′ be successive real zeros of n(s, χ) and σ be the sign of n′(u, χ). Restrict s = iw to be on the strip V(u, u′). σ/n(iw, χ) is positive definite in w.

### (8.2) The Ramanujan tau Dirichlet function, r(s).

See T. M. Apostol [4], G. H. Hardy [17], Eric W. Weisstein [42].

Let τ(n) be generated by v·$(\Pi_{n \geq 1}(1 - v^n))^{24} = \sum_{n \geq 1} \tau(n) \cdot v^n$.

**Definition of the Ramanujan tau Dirichlet function, r(s).**

r(s), for Re(s) > 7, is given by

$$r(s) := \sum_{n \geq 1} \tau(n) \cdot n^{-s} = \Pi_{p \text{ prime} > 0} (1 - \tau(p) \cdot p^{-s} + p^{11 - 2s})^{-1} \text{ (Euler product)}.$$

r(s) extends to an entire function of s on C. Define q(s) := $(2\pi)^{-s} \Gamma(s) r(s)$. The functional equation is

$$q(12 - s) = q(s).$$

**Ramanujan conjecture, RC.** Each nonreal zero of the Ramanujan function, r(s), has real part 6.

See J. Keiper [22].

**Simple zeros conjecture, SZCR.** The nonreal zeros of r(s) are all simple.

The q associated with the Ramanujan function, r(s), is used to define n(s) := sin(πs/2)q(6 + s). This n(s) is entire and odd.

The next conjecture restricted to 0 < x < 1 implies RC and SZCR.

**Conjecture R.** n(s) is a groove function in s on the complex plane.

**Main conjecture for r(s).** Let k be an integer. Restrict s = iw to the vertical strip V(2k, 2(k + 1)). Take σ to be the sign of n′(2k). Then σ/n(iw) is positive definite in w.

**Appendix**

**Kronecker's approximation theorem and Dirichlet L-functions, L(s, χ).**

**Introduction**

Our aim is to prove Corollary 3 below. (See Lemma 8.1 of Part I, §8). That establishes the existence of characters χ as follows. Say Re(s) > 1. L(s, χ) is not, even in the extended sense, either a groove or a ridge function (see Part I, §7 and §8).

**Definitions**

**Definitions of: invasive, discretely invasive (d-invasive) and $P_N(a)$.**
A vector in $R^N$ is said to be invasive (respectively: discretely invasive or d-invasive) when its coordinates (respectively: together with $2\pi$) are independent over the integers. An infinite sequence a of reals a(k) is invasive (respectively d-invasive) when each of its initial sections $P_N(a) := (a(1), a(2), ... a(N))$ is invasive (respectively d-invasive).

**Definition and invasiveness of ω.** The sequence ω of logarithms of the successive positive primes is invasive by reason of the uniqueness of the decomposition independent of order of any positive integer into a product of positive prime factors. ω is also discretely invasive, since the positive even powers of $e^\pi$ are irrational, $e^\pi$ being transcendental.

**Definition of $S_w(a, J)$.** Let w = 1, 2 and J > 0. If a is invasive but not d-invasive, take $S_w(a, J)$ to be the interval $(J, \infty)$ of reals greater than J. When a is d-invasive take $S_1(a, J) := (J, \infty)$ and $S_2(a, J)$ to be the set of integers greater than J.

**Definitions of T and $d_N$.** Let T be the abelian group of the reals mod($2\pi$) under addition. Take the metric $d_N$ on the N-D torus $T^N$ to be $d_N(u, v) := \max\{|s(k)| : 1 \le k \le N\}$, with $s(k) \equiv (u(k) - v(k)) \bmod(2\pi)$ and $-\pi < s(k) \le \pi$.

Let z be an infinite sequence of complex numbers z(k). Henceforth assume $|z(k)| < 1$ and $\sum_{k \ge 1} |z(k)|$ converges. Fix z.

**Definitions of $g_k(\theta)$, $G_N(y)$ and G(y).** $g_k(\theta) := |1 - z(k)e^{i\theta}|$. Let y be in $T^N$. Set $G_N(y) := \Pi_{1 \le k \le N} \, g_k(y(k))$. Say y is in $T^\infty$. Take $G(y) := \Pi_{k \ge 1} \, g_k(y(k))$.

Note that for real θ: $0 < 1 - |z(k)| \le g_k(\theta) \le 1 + |z(k)|$. $\Pi_{k \ge 1} (1 - |z(k)|)$ converges and is positive. $\Pi_{k \ge 1} (1 + |z(k)|)$ converges.

The sequence of functions $G_N$, with $N \geq 1$, is equi-uniformly continuous: Given $\varepsilon > 0$, there exists $\delta > 0$ such that for any $N \geq 1$ and u, v on the compact N-dimensional torus $T^N$ with $d_N(u, v) < \delta$, one has $|G_N(u) - G_N(v)| < \varepsilon$.

As N becomes infinite, $G_N(P_N(y))$ converges to $G(y)$ uniformly in y on $T^\infty$.

Say y, y' are on $T^\infty$ and $N \geq 1$.

**Definitions of $V_N(y)$ and $D_N(y', y)$.**
$V_N(y) := G(y) - G_N(P_N(y))$. $D_N(y', y) := G_N(P_N(y')) - G_N(P_N(y))$.

**Lemma 1** *Say $\varepsilon > 0$. There exist a positive integer $N(\varepsilon)$ and a $\delta(\varepsilon) > 0$ such that if $N \geq N(\varepsilon)$, x, y are in $T^\infty$ and $d_N(P_N(x), P_N(y)) < \delta(\varepsilon)$, then $|G(x) - G(y)| < \varepsilon$.*
.
**Proof** Take $N(\varepsilon)$ such that for any $N \geq N(\varepsilon)$ and x in $T^\infty$: $|V_N(x)| < \varepsilon/3$. There exists a $\delta(\varepsilon) > 0$ such that for any $N \geq N(\varepsilon)$ and u, v in $T^N$ with $d_N(u, v) < \delta(\varepsilon)$, one has $|G_N(u) - G_N(v)| < \varepsilon/3$.

Set $r(1) := V_N(x)$, $r(2) := D_N(x, y)$ and $r(3) := -V_N(y)$. Then $G(x) - G(y) = r(1) + r(2) + r(3)$. Also $|r(k)| < \varepsilon/3$. Hence $|G(x) - G(y)| < \varepsilon$.

Use will be made of the Kronecker approximation theorem (see T. M. Apostol [4]) in the following form.

**Kronecker approximation theorem.** *If $\alpha$ in $T^N$ is invasive (respectively: d-invasive) and $J > 0$, then the $t\alpha$, arising from the real (respectively: integral) t with $t > J$, are dense in $T^N$.*

**Corollary 1** *Suppose $\alpha$ in $T^\infty$ is invasive. Say $\varepsilon$, J are positive. There exist a positive integer N' and a $\delta > 0$ such that, for any sequence y in $T^\infty$ and integer $N \geq N'$, the following hold.*
*(1) There exists t in $S_w(\alpha, J)$ with $d_N(t \cdot P_N(\alpha), P_N(y)) < \delta$.*
*(2) If t is as in (1), then $|G(t\alpha) - G(y)| < \varepsilon$.*

**Lemma 2** *Suppose $\alpha$ in $T^\infty$ is invasive, $\sigma = \pm 1$ and $J > 0$. Then there exists a t in $S_w(\alpha, J)$ with $\Delta' := G(t\alpha) - G(0)$ nonzero and of the same sign as $\sigma$.*

**Proof** Say $\theta$ is real and for some k': $\Delta := g_{k'} \cdot (\theta) - g_{k'} \cdot (0)$ is nonzero and has the same sign as $\sigma$. Let $\delta(k', k) := 0$, except for $\delta(k', k') := 1$. $G(\theta \delta(k', k)) - G(0) = \Phi \Delta$, with $\Phi := \Pi_{k \geq 1, k \neq k'} g_k(0)$. $\Phi$ is positive, since $\Phi > \Pi_{k \geq 1} (1 - |z(k)|) > 0$.

Corollary 1 assures *a fortiori* that for any positive $\varepsilon$, J there exists t in $S_w(\alpha, J)$ for which $|\Delta''| < \varepsilon$, with $\Delta'' := G(t\alpha) - G(\theta \delta(k', k))$. Now $\Delta' = \Delta'' + \Phi \Delta$. Take $\varepsilon$

less than $\Phi|\Delta|$.

**Corollary 2** *Say α in $T^\infty$ is invasive. Assume at least one of (1), (2) which follow holds.*
*(1) There is a k' with z(k') nonreal.*
*(2) There exist k', k'' with z(k') < 0 and z(k'') > 0.*
*Then for any J > 0 there exist t, t' in $S_w(α, J)$ with G(tα) > G(0) > G(t'α).*

**Corollary 3** *Assume that either (1) or (2) below holds.*
*(1) χ(n) is nonreal, for some n.*
*(2) There exist primes p, p' such that χ(p) = 1 and χ(p') = -1.*

*Then for any x > 1 and T > 0, there exist real numbers (integers) t, t' > T, such that:*

$$|L(x + it, χ))| < |L(x, χ)| < |L(x + it', χ)|.$$

**Proof** Say $s = x + iθ$ with $x > 1$ and θ real. The Euler factorization for $L(s, χ)$ (see Part I, §8, and T. M. Apostol [3]) yields $|L(s, χ)| = 1/G(-θω)$, with G arising from the sequence z with kth term $χ(p_k)·(p_k)^{-x}$, where $p_k$ is the kth positive prime. ω is d-invasive. Corollary 2 yields Corollary 3.

**Acknowledgements**

An initial part of this work was done during a stay at the *Institut des Hautes Études Scientifiques*, France, as Visiting Professor. The author expresses his gratitude to Professor Maxim Kontsevich for a pleasurable collaboration which benefited from his deep mathematical insight and encouragement.

The author thanks Cécile Cheikhchoukh for word processing a manuscript presenting that initial phase of research. Appreciation is also extended to the IHES for its hospitality.